\documentclass[12pt]{amsart}
\usepackage[utf8]{inputenc}
\usepackage{amsmath,amsthm,amscd,amssymb,enumerate,enumitem}
\usepackage[margin=1in]{geometry}
\usepackage{color}
\usepackage{fancyhdr}
\usepackage{amsfonts}
\usepackage{mathtools}
\usepackage{etoolbox}
\usepackage{amsfonts}
\usepackage{cases}
\usepackage{tikz-cd}
\usepackage{faktor} 
\usepackage[backend=biber,style=numeric,sorting=ynt]{biblatex}
\usepackage[colorlinks=true, linkcolor=blue, citecolor=blue, urlcolor=blue]{hyperref}

\usepackage{esint}
\numberwithin{equation}{section}
\newtheorem{Theorem}{Theorem}[section]
\newtheorem{Proposition}[Theorem]{Proposition}
\newtheorem{Lemma}[Theorem]{Lemma}
\newtheorem{Corollary}[Theorem]{Corollary}
\newtheorem{Claim}[Theorem]{Claim}
\theoremstyle{definition}
\newtheorem{Definition}[Theorem]{Definition}

\newtheorem{Remark}[Theorem]{Remark}

\title[Diff. Finiteness for Non-Collapsed Ricci Flows on K\"ahler Surfaces]{Diffeomorphism Finiteness for Non-Collapsed Solutions to the Ricci Flow on K\"ahler Surfaces}
\author[Fluck]{Harry Fluck}
\address{Department of Mathematics, University of Illinois Urbana-Champaign \\ Urbana, IL 61801, USA}
\email{hpf5@illinois.edu}
\date{}

\begin{document}

\begin{abstract}
    In this article, we prove a diffeomorphism finiteness result for non-collapsed solutions to the Ricci flow on K\"ahler surfaces. 
     Our proof is via a neck decomposition type theorem for $P^\ast$-parabolic balls below the type-I scale.
    As a consequence of our decomposition, we recover certain sharp Minkowski estimates for the almost singular set originally due to Fang and Li. 
\end{abstract}

\maketitle

\section{Introduction}

The Ricci flow was first introduced by Richard Hamilton \cite{3manPRc} in 1982 and has since proven to be a powerful tool in geometry and topology.
 In general, solutions to this flow will encounter finite time singularities, and hence understanding the possible structures of blowup models is of great interest.   In dimensions $2$ and $3$, singularity models are smooth and have been completely classified \cite{brendle2020,Brendle_2021,RF_surfaces,RF_surfaces_2}. In dimension $4$ and higher, singularity models are potentially singular and at present any complete classification is out of reach. As a weaker question, it is natural to ask whether the small scale topology of an arbitrary solution to the flow may at least be bounded.  In dimensions $5$ and higher, this is already known to be false even in the elliptic setting \cite{anderson_diff}. In dimension $4$, it was shown in \cite{regEinstein} that the set of uniformly non-collapsed $4$-manifolds with bounded Ricci curvature and bounded diameter have bounded diffeomorphism type, resolving a conjecture of Anderson. Using these ideas, it was shown in \cite{Bubble_Decomp} that the set of closed $4$-dimensional Ricci flows with bounded scalar curvature and bounded diameter have bounded diffeomorphism type. The crucial simplification is that the scalar curvature bounds imply $2$-sided heat kernel bounds and distance distortion estimates, allowing one to work on a fixed time slice. This approach cannot be extended to general $4$-dimensional Ricci flows, where the question remains open. Recently, Richard Bamler introduced a compactness and partial regularity theory that allows for the study of singular solutions to the Ricci flow in great generality. This theory exhibits remarkable similarities with that for Gromov-Hausdorff limits of Einstein metrics, and so it is natural to ask whether these techniques may be combined with those of \cite{regEinstein} to bound the small scale topology of $4$-dimensional Ricci flows.

 In this paper, we will show that below the type-I scale, away from a set of times of small measure, the diffeomorphism type of \textit{any} non-collapsed solution to the Ricci flow on a closed K\"ahler surface is bounded. In particular, our main result is the following.

 \begin{Theorem}[Finite Diffeomorphism Type]\label{thm:diffeo}
    Let $\epsilon>0$ and $Y<\infty $. There exists $\Lambda(Y),C(Y,\epsilon)<\infty$ such that the following holds. Let $(M^4,(g_t)_{t\in [-10,10]})$ be a closed K\"ahler Ricci flow and  $x_0\in M $ with $\mathcal{N}_{x_0,0}(1)\geq -Y$. There exists $K\subseteq [-1,1]$ with $|K|<\epsilon$ such that for all $t\in [-1,1]\backslash K$,  there is an open set $\Omega_t\subseteq M $ with $$P_t^\ast(x_0,0;1) \subseteq \Omega_t \subseteq P^\ast_t(x_0,0;\Lambda,2,2),$$ 
    which has one of $C(Y,\epsilon)$ diffeomorphism types. 
\end{Theorem}

The proof of Theorem \ref{thm:diffeo} is inspired by the arguments in \cite[Section 8.4]{regEinstein}. In particular, we produce a neck-body decomposition, which states that a $P^{\ast}$-parabolic ball below the type-I scale may be decomposed into a collection of regions of two types: annular regions, on which the flow is close to an annulus in the static flow modeled on a flat cone of the form $C(\mathbb{S}^3/\Gamma)$, and body regions, on which the curvature scale is uniformly controlled by the diameter. In particular, we will prove the following (c.f. \cite[Theorem 8.64]{regEinstein} and \cite[Theorem 3.2]{Bubble_Decomp}).

\begin{Theorem}[Annular Decomposition] \label{thm:neck_decomp}
    Let $Y<\infty$. Then the following is true for  $\Lambda\geq \underline{\Lambda}(Y)$, $\gamma\leq \overline{\gamma}(\Lambda,Y)$, and $\epsilon \leq \overline{\epsilon}(Y,\gamma,\Lambda)$. Let $(M^4, (g_t)_{t\in [-10,10]})$ be a closed K\"ahler Ricci flow and consider $x_0\in M$ with $\mathcal{N}_{(x_0,0)}(1)\geq -Y$. There exists a covering $$ P^\ast(x_0,0;1)\subseteq \bigcup_{i\in N}\mathcal{A}^i\cup \bigcup_{j\in \hat{N}}\mathcal{B}^j\subseteq P^\ast(x_0,0;\Lambda,2,2)$$
    such that 
    \begin{enumerate}
        \item for each $i\in N$, $\mathcal{A}^i=P^\ast(x_i,t_i;\Lambda r_i,2r_i^2,2r_i^2)\backslash \bigcup_{c\in\mathcal{C}^i}P^\ast(x_c,t_c;r_c)$ is an $(m,\epsilon,\gamma,\Lambda)$-annular region for some $m\leq e^{Y}$,
        \item  for each $j\in\hat{N}$, we have  $\mathcal{B}^j=\bigcup_{e\in E^j}P^\ast(x_e,t_e;r_e)$  where for each $e\in E^j$ we have $r_{\operatorname{Rm}}\vert_{P^\ast(x_e,t_e;r_e)}\geq c(Y,\epsilon)r_e$,
    \end{enumerate}
     and we have the content estimate $$ \sum_{i\in N} r_i^2+ \sum_{j\in \hat{N}}\sum_{e\in E^j}r_e^2\leq C(Y,\epsilon).$$
\end{Theorem}

\begin{Remark}
   Roughly speaking, $\mathcal{A}^i$ is an $(m,\epsilon,\gamma,\Lambda)$-annular region if the flow is close to that of a static flow modeled on the annulus of a flat cone of the form $C(\mathbb{S}^3/\Gamma)$ on the entire region. For the precise definition of an $(m,\epsilon,\gamma,\Lambda)$-annular region, see Definition \ref{def:annular}. Note that we cannot bound the actual \textit{number} of annular/body regions in our decomposition. This is due to the fact that almost singular points are only isolated on a \textit{fixed} time slice, but will form a collection of approximate temporal lines in spacetime. We also remark that a more general neck decomposition type theorem for non-collapsed $4$-dimensional Riemannian Ricci flows was given in \cite[Theorem 6.9]{fang2026singularsetsnoncollapsedricci}.
\end{Remark}

Based on our construction, it will follow from standard theory that the diffeomorphism type of the time slices of the body regions is bounded (c.f. \cite[Theorem 8.4]{regEinstein}), whilst that of the annular regions is controlled by the uniform non-collapsing estimate. Modulo some technical lemmas, the proof of Theorem \ref{thm:diffeo} follows by combining these two facts with Theorem \ref{thm:neck_decomp} in the same way as \cite[Proof of Theorem 1.12]{regEinstein}.

The proof of Theorem \ref{thm:neck_decomp} relies on a new observation that rules out the appearance of cylindrical necks below the type-I scale for non-collapsed Ricci flows on K\"ahler surfaces (c.f. Lemma \ref{lem:weaksplittingepsreg}). This implies that, below the type-I scale, the almost singular set is approximately isolated in space, which allows for a construction similar to that of \cite{regEinstein, Bubble_Decomp}. That being said, due to the lack of assumed curvature bounds in Theorem \ref{thm:diffeo}, we cannot work on a fixed time slice, as the Riemannian distance function is not well controlled. As a result, we are forced to work in spacetime, using the natural topology of the flow, and restrict to a general time slice only at the very end of our proof, which complicates things significantly. In particular, we will follow an argument similar to that in \cite{L2CurvEin}, where one employs techniques from rectifiable Reifenberg theory.
More specifically, we introduce the notion of $(c), (d)$, and $(e)$-ball (c.f. Definition \ref{def:balls} and the proceeding discussion), which are types of $P^\ast$-parabolic balls differentiated by the behaviour of their entropy pinched points. Our first step is then to prove an initial covering lemma for balls of both $(c)$ and $(d)$ type (Lemmas \ref{lem:coarse_c} and  \ref{lem:coarse_d}). Roughly speaking, these lemmas state that a $(c)$-ball has the structure of an $(m,\epsilon,\gamma,\Lambda)$-annular region, whilst a $(d)$-ball may be decomposed into $(e)$-balls, modulo a collection of balls of very small content. The main technical difficulty, which does not appear in the elliptic setting, is the presence of so called quasi-static regions in these covers, where a ball is close to the static flow modeled on a flat cone only up to some fixed time. In light of Theorem \ref{thm:diffeo}, the diffeomorphism type of the time slices of such regions cannot be controlled past this time. To circumvent this problem, we introduce an additional ball type, known as an $(f)$-ball (c.f. Definition \ref{def:balls}), on which this issue occurs. We then prove a covering lemma for such balls (c.f. Lemma \ref{lem:f_ball}), which states that they may be decomposed into  $(c), (d),$ and $(e)$-balls, whilst retaining control on the $2$-content. Combining this fact with our initial coverings, we then use an induction argument to construct refined $(c)$ and $(d)$-ball coverings (c.f. Lemmas \ref{lem:ref_c} and \ref{lem:ref_d}), which contain no quasi-static regions. The key observation in this step is that since the almost singular points form an aproximate temporal line, both the $(f)$-ball covering, and the refined $(c)$-ball covering, can be reduced to a simple packing argument on the real line (c.f. Lemma \ref{lem:real_analysis}). Using these improved covering lemmas, we then prove Theorem \ref{thm:neck_decomp} via a final induction argument involving the $(e)$-balls. This component of the argument is largely unchanged from \cite[Proof of Theorem 7.1]{L2CurvEin}.  We  expect these new techniques will be useful, not only in the study of rectifiable Reifenberg theory for the Ricci flow in more general settings, but also for the corresponding theory for other parabolic equations.


As an application of Theorem \ref{thm:neck_decomp}, we recover certain sharp Minkowski estimates for the almost singular set, originally due to Fang and Li (more specifically, we recover \cite[Theorem 1.5]{fang2026singularsetsnoncollapsedricci} in the K\"ahler case).

\begin{Theorem}[Sharp Minkowski Estimates]\label{thm:minkowski}
Let $Y<\infty$. Let $(M^4, (g_t)_{t\in [-10,10]})$ be a closed K\"ahler Ricci flow, and let $x_0\in M $ be such that $\mathcal{N}_{(x_0,0)}(1)\geq -Y$. Then for all $\sigma>0$, we have
$$|\{(x,t)\in P^\ast(x_0,0;1)\,\vert\, r_{\operatorname{Rm}}(x,t)<\sigma \}|\leq C(Y)\sigma ^4.$$
\end{Theorem}

Before concluding this section, we remark that the K\"ahler assumption in Theorem \ref{thm:diffeo} is expected to be strictly necessary. For Riemannian Ricci flow, we anticipate that Lemma \ref{lem:weaksplittingepsreg} will not hold (informally, this can be seen by considering, for instance, a very small scale degenerate neck pinch), and so an additional region type, one on which the flow is close to a generalized cylinder, would be present in Theorem \ref{thm:neck_decomp} (see \cite[Theorem 6.9]{fang2026singularsetsnoncollapsedricci}). In particular, almost singular points are no longer isolated in space.

The remainder of our paper is structured as follows. In Section \ref{sec:preliminaries}, we discuss some technical preliminaries that will be needed throughout this article. In Section \ref{sec:annular}, we introduce the notion of an $(m,\epsilon,\gamma,\Lambda)$-annular region, and prove some important properties pertaining to their structure. In Section \ref{sec:ball}, we introduce the various types of balls that will be used in our covering arguments. In Section \ref{sec:toy}, we prove the packing estimate on the real line that will play a role in the proofs of Lemma \ref{lem:f_ball} and Lemma \ref{lem:ref_c}. In Section \ref{sec:neck_decomp}, we prove Theorem \ref{thm:neck_decomp}. In Section \ref{sec:main_thm}, we prove Theorem \ref{thm:diffeo} and Theorem \ref{thm:minkowski}.

\subsection*{Acknowledgments} The author is very grateful to Max Hallgren, for many insightful discussions related to this paper. The author would also like to thank both Xiaodong Cao and Eric Chen, for their support and encouragement.

\subsection*{AI-Usage Statement}The mathematical ideas, arguments, and text, in this article were produced without any assistance from large language models. Claude and ChatGPT were used in the final stages of preparation in order to conduct literature searches, check for grammatical errors and notational inconsistencies, and to verify the validity of some technical components of the argument. In this stage, some minor technical amendments suggested by Claude were employed. The author assumes full responsibility for the mathematical content of this work.

\section{Preliminaries}\label{sec:preliminaries}

In this section, we fix notation and introduce some important results and definitions that will be used throughout this article.

\subsection{Notation}

Throughout the remainder of the paper, we adhere to the following notational conventions. The notation $A<\infty$ means that $A$ is a large positive constant, while $\epsilon>0$ means that $\epsilon$ is a small positive constant. We let $\Psi(a_1 ,...,a_k|b_1,...,b_{\ell})$ denote a quantity depending on parameters $a_1,...,a_k,b_1,...,b_{\ell}$, which satisfies
$$\lim_{(a_1,...,a_k)\to (0,...,0)} \Psi(a_1,...,a_k|b_1,...,b_{\ell})=0$$
for any fixed $b_1,...,b_{\ell}$. When we write $Q\leq C(a_1,...,a_k)$ we mean that a given quantity $Q$ may be bounded from above by some constant depending on the parameters $a_1,...,a_k$. The expression $Q\geq c(a_1,...a_k)$ is defined similarly. Also, if we say that a proposition $P(\epsilon)$ depending on a parameter $\epsilon$ holds if $\epsilon \leq \overline{\epsilon}(b_1,...,b_{\ell})$, this means there exists a constant $\overline{\epsilon}$ depending on parameters $b_1,...,b_{\ell}$ such that $P(\epsilon)$ holds whenever $\epsilon \in (0,\overline{\epsilon}]$. The notation $E \geq \underline{E}(b_1,...,b_{\ell})$ is defined analogously. The precise value of the constants present in our results will never be of relevance. As a result, when estimating, we will allow constants to change from one line to the next, without introducing new notation. We also omit dependence on the dimension $n$ throughout.

\subsection{Setup and Basic Concepts}

Our central object of study will be closed Ricci flows $(M^n,(g_t)_{t\in I})$, where $M^n$ is a closed $n$-dimensional manifold, $I\subseteq \mathbb{R}$ is an interval, and $(g_t)_{t\in I}$ is a smooth family of Riemannian metrics on $M$ which satisfy $$\partial_t g_t=-2Rc_{g_t}.$$ We will be most interested in the case where $M$ is K\"ahler and of complex dimension $2$.

We let $\mathfrak{t}: M\times I \to I$ denote the projection onto the time factor. We let $B(x,t;r)$ denote the usual geodesic ball of radius $r>0$ centered at $x\in M$, with respect to the metric $g_t$. We define the curvature scale of $(x,t) \in M\times I$ to be 
$$r_{\operatorname{Rm}}(x,t)\coloneqq \sup \left\{ r>0\,\vert \, |Rm|(x',t')\leq r^{-2} \text{ for all } (x',t') \in B(x,t;r) \times ([t-r^2,t+r^2]\cap I)  \right\}.$$ Given  $Y\subseteq M$, we let $|Y|_{g_t}$ and $\operatorname{diam}_{g_t}(Y)$ denote the volume and diameter of $Y$  with respect to the metric $g_t$, respectively. Given a space time subset $\widetilde{Y}\subseteq M\times I$, we let $$|\widetilde{Y}| \coloneqq \int_I |\widetilde{Y}\cap (M\times \{t\})|_{g_t}\, dt $$ denote the parabolic volume.

\subsection{Conjugate Heat Kernels}

We let
$$\Box \coloneqq  \partial_t-\Delta_{g_t}, \qquad \Box^{\ast}\coloneqq -\partial_t-\Delta_{g_t}+R_{g_t}$$
denote the heat operator and its formal adjoint, respectively. Given $(x_{0},t_{0})\in M \times I$ the \textit{conjugate heat kernel} based at $(x_0,t_0)$ is the function $K(x_{0},t_{0};y,s)\in C^{\infty}(M\times (-\infty,t_{0})\cap I)$ defined by $$ \begin{cases}
        \lim_{s\nearrow t_0}  K(x_0,t_0;y,s)=\delta_{x_0} \\ \Box^{\ast}_{y,s}K(x_0,t_0;y,s)=0
    \end{cases}.$$
The weighted measures $\nu_{x_{0},t_{0};s}\coloneqq K(x_{0},t_{0};\cdot ;s)dg_{s}$ are probability measures, known as the \textit{conjugate heat kernel measures} based at $(x_0,t_0)$. It is often convenient to write $K(x_0,t_0;y,t)=(4\pi\tau)^{-n/2}e^{-f}$, where $\tau\coloneqq t_{0}-t$, and $f\in C^{\infty}(M)$. From this, there is an associated quantity
\begin{align} \label{eq:nashentropy}
    \mathcal{N}_{x_0,t_0}(\tau)\coloneqq \int_{M} f(y,t_0-\tau)\,d\nu_{x_0,t_0;t_0-\tau}(y)-\frac{n}{2},
\end{align}
known as the \textit{pointed Nash entropy} based at $(x_0,t_0)$. It is well known that for each $(x_0,t_0)\in M\times I$, the map $\tau \xrightarrow[]{} \mathcal{N}_{x_0,t_0}(\tau)$ is non-increasing \cite[Proposition $5.2$]{Bam1} . We will often assume a bound of the form $\mathcal{N}_{(x_0,t_0)}(r^{2})\geq -Y$, which can be interpreted as a volume-noncollapsing condition near $(x_0,t_0)$ at scale $r$ (c.f. \cite[Theorem 6.1 and Theorem 8.1]{Bam1}). Due to a variation bound, after adjusting constants, we may always pass pointwise Nash entropy bounds to local bounds (c.f. \cite[Corollary 5.11]{Bam1}).

We refer the reader to \cite[Sections 2 and 3]{Bam1} for the definition of the Wasserstein distance $d_{W_1}$ and related notions, and their properties on a Ricci flow. We let $d_{W_1}^{g_t}$ denote the Wassserstein distance with respect to the metric $g_t$. Given $(x_0,t_0), (x_1,t_1)\in M\times I$, the map
  \begin{align*}
      t\to d_{W_1}^{g_t}(\nu_{x_0,t_0;t},\nu_{x_1,t_1;t})
  \end{align*}
  is non-decreasing \cite[Lemma $2.7$]{Bam1}. For this reason, the following is a well-behaved replacement for the usual notion of parabolic balls.
  \begin{Definition} \label{def:Pstarparabolic}
      Let $(x_0,t_0)\in M\times I$. For $A,T^+,T^->0$, we define  $$P^{\ast}(x_0,t_0;A,T^-,T^+)=\{ (x,t)\in M\times [t_0-T^-,t_0+T^+]\, \arrowvert \,d_{W_1}^{g_{t_{0}-T^-}}(\nu_{x_0,t_0;t_0-T^-},\nu_{x,t;t_0-T^-})<A\}.$$ 
      We refer to $P^{\ast}(x_0,t_0;r)\coloneqq P^{*}(x_0,t_0;r,r^2,r^2)$ as the $P^*$-parabolic ball of radius $r>0$ around $(x_0,t_0)$. We also set $P^\ast_t(x_0,t_0;r)\coloneqq P^\ast(x_0,t_0;r)\cap (M\times \{t\})$, and given $J\subseteq I$, we set $P^\ast_J(x_0,t_0;r)=P^\ast(x_0,t_0;r)\cap (M\times J)$.
      
  \end{Definition}
 Given a spacetime subset $Y\subseteq M\times I$, we set $$\operatorname{diam}_{P^\ast}(Y)=\inf \{ r>0\, \vert\, \exists (x,t)\in M\times I\, \text{ s.t. } P^\ast(x,t;r)\supseteq Y\}.$$ We define the $\beta$-tubular neighborhood around $Y$ with respect to the natural topology to be $$T_{\beta}(Y)=\bigcup_{(x,t)\in Y}P^\ast(x,t;\beta).$$  Note that $P^\ast$-parabolic balls satisfy a triangle type inequality with slightly worse constants (c.f. \cite[Proposition 9.4]{Bam1}). Due to a Vitali covering lemma \cite[Lemma 17.35]{Bam3}, and the following parabolic volume doubling esimate, we may apply standard covering arguments to $P^\ast$-parabolic balls (see \cite[Lemma 6.1]{fluck2025epsilonregularitytheoremnoncollapsedricci}, or \cite[Theorem 9.8]{Bam1} for the upper bound).
\begin{Lemma} \label{lem:non_collapsing}
	For every $Y<\infty$, there exists a constant $C(Y)<\infty$ such that the following holds.    Let $(M^n,(g_t)_{t\in I})$ be a Ricci flow on a closed manifold, and suppose $(x_{0},t_{0})\in M\times I$, $r>0$ satisfy $[t_0-2r^2,t_0]\subseteq I$ and $\mathcal{N}_{(x_{0},t_{0})}(r^{2})\geq -Y$. Then for all $s\in [0,r]$ we have 
	\begin{align*}
		C(Y)^{-1}s^{n+2}\leq  |P^{\ast}(x_{0},t_{0};s)|\leq C(Y)s^{n+2}.
	\end{align*}
\end{Lemma}

 The conjugate heat kernel measure based at $(x_0,t_0)$ may be viewed as a probability distribution for the location of $x_0$ at previous times. The following notion can be thought of as a choice of ``mean" for $\nu_{x_0,t_0;t}$. 
\begin{Definition}
	Let $(x_0,t_0)\in M\times I$, and set $H_n\coloneqq \frac{(n-1)\pi^2}{2}+4$. A point $(z,t)\in M\times I\cap (-\infty,t_0]$ is called an $H_n$-center of $(x_0,t_0)$ if 
	$$\int_M d_{g_t}^2(z,y)d\nu_{x_0,t_0;t}(y) \leq H_n(t_0-t).$$
\end{Definition}
It follows from a concentration bound that $H_n$-centers always exist at each time $t\leq t_0$  \cite[Proposition $3.12$]{Bam1}, and that any two $H_n$-centers $(z,t),(z',t)$ satisfy $d_{g_t}(z,z')\leq 2\sqrt{H_n(t_0-t)}$.

\subsection{Integral Almost Properties}

 Next we recall the integral-almost properties introduced in \cite{Bam3} and \cite{HJ}. These measure how close the flow is to a particular model solution near a given point. 
 

Throughout this subsection, we let $(M,(g_t)_{t\in I})$ denote a compact $n$-dimensional Ricci flow. 

\begin{Definition} \label{def:selfsimilar}
     Let $(x_0,t_0)\in M \times I$, $\epsilon,r>0$, and write $W\coloneqq \mathcal{N}_{x_0,t_0}(r^2)$. We say that $(x_0,t_0)$ is $(\epsilon ,r)$-selfsimilar if $[t_0-\epsilon^{-1}r^2,t_0]\subseteq  I $ and 
         \begin{enumerate}
             \item \label{iden:selfsimilar1} $\int_{t_0-\epsilon^{-1}r^{2}}^{t_0-\epsilon r^{2}}\int_M \tau |Rc+\nabla^{2}f-\frac{1}{2\tau }g|^{2}\,d\nu_{x_0,t_0;t}\,dt<\epsilon$,
             \item \label{iden:selfsimilar2} $\sup _{t\in [t_0-\epsilon^{-1}r^{2},t_0-\epsilon r^{2}]}\int_M |\tau(R+2\Delta f-|\nabla f|^{2})+f-n-W|\,d\nu_{x_0,t_0;t}<\epsilon,$
             \item \label{iden:selfsimilar3} $\inf_{M\times [t_0-\epsilon^{-1}r^{2},t_0-\epsilon r^{2}]}R(x,t)\geq -\epsilon r^{-2}$
         \end{enumerate}
   where $f$ is given by $d\nu_{x_0,t_0;t}=(4\pi \tau)^{-\frac{n}{2}}e^{-f}dg_t$.
\end{Definition}

\begin{Definition} \label{def : strong_splitting}
    Let $(x_0,t_0)\in M \times I$ and $\epsilon,r>0$. We say that $(x_0,t_0)$ is strongly $(k,\epsilon,r)$-split if $[t_0 - \epsilon^{-1}r^2,t_0] \subseteq I$ and there exist $y_1,...,y_k \in C^{\infty}(M \times [t_0 - \epsilon^{-1}r^2,t_0-\epsilon r^2])$ such that the following hold for $1\leq i,j\leq k$:
    \begin{enumerate}
        \item $\Box y_i=0, \label{iden : strong_splitting1}$ 
        \item $r^{-2} \int_{t_0-\epsilon^{-1}r^2}^{t_0-\epsilon r^2} \int_M |\langle \nabla y_i, \nabla y_j\rangle - \delta_{ij}| d\nu_{x_0,t_0;t}dt\leq \epsilon$ \label{iden : strong_splitting2},
        \item $\int_M y_i \, d\nu_{x_0,t_0;t}=0 \text{ for each }t\in [t_0-\epsilon^{-1} r^2,t_0-\epsilon r^2]$ \label{iden : strong_splitting3}.
    \end{enumerate}
\end{Definition} 

\begin{Definition} \label{def : static }
    Let $(x_0,t_0)\in M\times I$ and  $\epsilon,r>0$. We say that $(x_0,t_0)$ is $(\epsilon,r)$-static, if $[t_0-\epsilon^{-1}r^2,t_0]\subseteq I$ and the following holds 
    \begin{enumerate}
        \item $r^{2} \int_{t_0-\epsilon^{-1}r^{2}}^{t_0-\epsilon r^{2}}\int_M|Rc|^2\,d\nu_{x_0,t_0;t}\,dt<\epsilon$ \label{iden : strong_static1},
        \item $\sup_{[t_0-\epsilon^{-1} r^2, t_0-\epsilon r^2]}\int _{M} R\, d\nu_{x_0,t_0;t}<\epsilon$ \label{iden : strong_static2},
        \item $\inf_{M\times [t_0-\epsilon^{-1} r^2, t_0-\epsilon r^2] }R(x,t)\geq -\epsilon r^{-2}$ \label{iden : strong_static3}.
    \end{enumerate}
\end{Definition}

\subsection{Limits of Non-Collapsed Flows}

In this section, we recall some aspects of the compactness and partial regularity theory of \cite{Bam1,Bam2,Bam3} that will be needed throughout this article. We only give the rough ideas behind the definitions, we refer the reader to \cite{Bam1,Bam2,Bam3} for a detailed account. In \cite{Bam2} a topology is constructed, known as the $\mathbb{F}$-topology, such that for \textit{any} sequence of closed  Ricci flows $(M^n_i,(g_{i,t})_{t\in I})$, with conjugate heat kernels $(\nu^i_{x_{0,i},t_{0,i};t})_{t\leq t_{0,i}}$, we have $$(M^n_i,(g_{i,t})_{t\in I}),(\nu^i_{x_{0,i},t_{0,i};t})_{t\leq t_{0,i}})\xrightarrow[i\to\infty]{\mathbb{F}}((\mathcal{X}_t)_{t\leq t_{0,\infty}},(\nu^\infty_{x_\infty;t})_{t\leq t_{0,\infty}})$$
along a subsequence. Here, the limiting object $((\mathcal{X}_t)_{t\leq t_{0,\infty}},(\nu^\infty_{x_\infty;t})_{t\leq t_{0,\infty}})$ is not a smooth Ricci flow, but is something known as a metric flow. This is a $1$-parameter family of metric spaces together with a family of time dependent probability measures that behave similarly to the conjugate heat kernel coupled to a smooth Ricci flow. The $\mathbb{F}$-convergence is then approximately equivalent to GW1-convergence on almost every time slice.  Additionally, we can fix a correspondence $\mathfrak{C}$, such that  
$$(M^n_i,(g_{i,t})_{t\in I}),(\nu^i_{x_{0,i},t_{0,i};t})_{t\leq t_{0,i}})\xrightarrow[i\to\infty]{\mathbb{F},\mathfrak{C}}((\mathcal{X}_t)_{t\leq t_{0,\infty}},(\nu^\infty_{x_\infty;t})_{t\leq t_{0,\infty}})$$
within $\mathfrak{C}$. Roughly speaking, $\mathfrak{C}$ is a family of large ambient metric spaces in which the sequence can be uniformally embedded, and hence the existence of $\mathfrak{C}$ will allow us to make sense of things such as convergence of points along our sequence. We will, however, mostly not need to work in this level of generality, and will only require the special case in which the flows are $4$-dimensional, and are either $(\epsilon_i,1)$-selfsimilar and $(k,\epsilon_i,1)$-split, or $(\epsilon_i,1)$-selfsimilar and $(\epsilon_i,1)$-static.

\begin{Theorem}[c.f. Theorems 7.6, 9.12, 9.31 in \cite{Bam2}] \label{thm:bamconvergence} Suppose $(M_i^4,(g_{i,t})_{t\in I_i})$ is a sequence of pointed Ricci flows and $x_i \in M_i$, $\delta_i \searrow 0$, are such that $(x_i,0)$ are  $(\delta_i,1)$-selfsimilar and $\mathcal{N}_{x_i,0}(1)\geq -Y$ for some $Y<\infty$. Then we can pass to a subsequence to obtain a metric soliton $(\mathcal{X},(\mu_t)_{t\in (-\infty,0]})$ along with a correspondence $\mathfrak{C}$ such that we have the following uniform $\mathbb{F}$-convergence on compact time intervals:
\begin{align}
    (M_i^n,(g_{i,t})_{t\in (-T_i,0]},(\nu_{x_i,0;t})_{t\in I_i})\xrightarrow[i\to \infty]{\mathbb{F},\mathfrak{C}} (\mathcal{X},(\mu_t)_{t\in (-T,0)}) \label{eq:Fconvergence}.
\end{align}
Moreover, there is a homeomorphism $\mathcal{X} \cong X\times (-\infty,0]$ (where $\mathcal{X}$ is equipped with the natural topology as in \cite[Section 3.6]{Bam2}) for some singular space $(X,d)$ (in the sense of \cite[Definition 2.15]{Bam3}) with regular set $(\mathcal{R}_X,g_X)$, and that the following hold under the identification $\mathcal{X} \cong X\times (-\infty,0)$:
\begin{enumerate}
    \item the regular set $\mathcal{R}$ of the metric flow (see \cite[Definition 9.11]{Bam2} is identified with $\mathcal{R}_X \times (-\infty,0)$ as Ricci flow spacetimes (see \cite[Definition 9.1]{Bam2}), where there exists $f\in C^{\infty}(\mathcal{R})$ such that (setting $\tau:=|t|$)
    $$Rc_{g_t}+\nabla^2 f = \frac{1}{2\tau}g, \qquad d\mu_t = (4\pi \tau)^{-\frac{n}{2}}e^{-f}dg_t$$
    on $\mathcal{R}$, and if $\partial_{\mathfrak{t}}$ is the timelike vector field of $\mathcal{R}$, then $\partial_{\mathfrak{t}}-\nabla f$ is identified with the standard vector field $\partial_t$ on the time factor for some $f\in C^{\infty}(\mathcal{R})$,
    \item there are time-preserving open embeddings $\psi_i:\mathcal{R}\supseteq U_i \to M_i$, where $(U_i)_{i\in \mathbb{N}}$ is a precompact open exhaustion of $\mathcal{R}$, such that
    $$\psi_i^{\ast}g_i \to g, \qquad \psi_i^{\ast}f_i' \to f, \qquad (\psi_i)_{\ast}\partial_t \to \partial_{\mathfrak{t}}$$
    in $C_{\operatorname{loc}}^{\infty}(\mathcal{R})$. 
\end{enumerate}
If in addition, $(x_i,0)$ is $(k,\delta_i,1)$-split for some $0 \leq k\leq 4$, then there is a splitting $\mathcal{X}\cong \mathcal{X}^\prime \times \mathbb{R}^k$ for some metric soliton $\mathcal{X}^\prime$ where $\mathbb{R}^k$ is the standard flow on $\mathbb{R}^k$. If $(x_i,0)$ is also $(\delta_i,1)$-static, then we have $\mathcal{X}\cong C(\mathbb{S}^3/\Gamma)$ where $C(\mathbb{S}^3/\Gamma)$ is the static flow modeled on a flat cone of the form $C(\mathbb{S}^3/\Gamma)$ for some $\Gamma \in SO(4)$ with $|\Gamma|\leq e^Y$.
\end{Theorem}
\section{Annular Regions}\label{sec:annular}

In this section, we introduce the notion of an $(m,\epsilon,\gamma,\Lambda)$-annular region, and prove some  facts regarding their structure that will be needed throughout this article. 

By an  $(m,\epsilon,\gamma,\Lambda)$-annular region, we mean a spacetime subset which is $\epsilon$-close to an annulus in a flat cone of the form $C(\mathbb{S}^3/\Gamma)$, on every time slice, potentially down to different scales. We will say that a region is  $(m,\epsilon,\gamma,\Lambda)$-quasi annular if this property holds away from a set of times of positive measure. In particular, we have the following.

\begin{Definition} \label{def:annular}
    Let $\epsilon,\gamma >0$, $\Lambda<\infty$, and $m\in \mathbb{N}\backslash \{1\}$. Let $(M^4,(g_t)_{t\in [-2\epsilon^{-1},10]} )$ be a closed Ricci flow, $x_0\in M$, $\mathcal{C}\subseteq P^\ast(x_0,0;2\Lambda ,2,2)$   be a countable subset, and $r:\mathcal{C}\to \mathbb{R}_{>0}$ be a function.  We say that  $\mathcal{A}\coloneqq P^\ast(x_0,0;2\Lambda ,2,2)\backslash \bigcup_{(x,t)\in\mathcal{C}} P^\ast(x,t;r_{(x,t)})$ is an $(m,\epsilon,\gamma,\Lambda)$-annular region if 
    \begin{enumerate}
    \item $\{P^\ast(x,t; \gamma^3\cdot r_{(x,t)})\}_{(x,t)\in \mathcal{C}}$ are pairwise disjoint,  \label{iden:annular1}
        \item $(x,t)$ is $(\epsilon,s)$-static for every $(x,t)\in \mathcal{C}$ and $s\in [r_{(x,t)},1]$, \label{iden:annular2}
        \item $|\mathcal{N}_{(x,t)}(\epsilon\cdot s^2)-\operatorname{log}(\frac{1}{m})|<\epsilon$ for every $(x,t)\in \mathcal{C}$ and $s\in [r_{(x,t)},\epsilon^{-1}]$,  \label{iden:annular3}
        \item $\bigcup_{(x,t)\in\mathcal{C}}(t-\gamma^5 \cdot r_{(x,t)}^2,t+\gamma^5\cdot r_{(x,t)}^2)\supseteq [\inf\mathfrak{t}(P^\ast(x_0,t_0;1)),\sup\mathfrak{t}(P^\ast(x_0,t_0;1))]$,  \label{iden:annular4}
    \end{enumerate}
where we set $r_{(x,t)}=r(x,t)$.    If \eqref{iden:annular1}-\eqref{iden:annular3} hold, but \eqref{iden:annular4} fails, we say that $\mathcal{A}$ is an $(m,\epsilon,\gamma,\Lambda)$-quasi annular region.
\end{Definition}

 We refer to the set $\mathcal{C}$ in Definition \ref{def:annular} as the center set. We begin with the following technical lemma, a consequence of which is that the center set of an $(m,\epsilon,\gamma,\Lambda)$-region approximately forms a temporal line.
 
 
 \begin{Proposition}\label{prop:time_dif}
 	For every $\gamma>0$, and $\Lambda,Y<\infty $, the following statement is true for $\epsilon \leq \overline{\epsilon}(Y,\gamma,\Lambda)$. Let $(M^4,(g_t)_{t\in I})$ be a closed Ricci flow, and $(y_1,s_1),(y_2,s_2)\in M \times I$ be $(\epsilon,r)$-static and $(\epsilon,r)$-selfsimilar with $\mathcal{N}_{(y_i,s_i)}(r^2)\geq -Y$ and $\mathcal{N}_{(y_i,s_i)}(\epsilon\cdot r^2)\leq \operatorname{log}(\frac{1}{2})+\epsilon $ for some $r>0$. Then if  $(y_2,s_2)\in  P^\ast(y_1,s_1;\gamma^3r)^c\cap  P^\ast(y_1,s_1;\Lambda r)$, we  must have $|s_1-s_2|\geq c(Y)\gamma^6r^2$. 
 \end{Proposition}
 \begin{proof}
 	By means of parabolic rescaling and time translation, we may assume that $r=1$ and $0=s_1\geq s_2$. Assume that the claim was false, then  there exists sequences $c_i,\epsilon_i\searrow 0$ and closed Ricci flows $(M^4_i,(g_{i,t})_{t\in [-\epsilon_i^{-1},0]})$ with points $(y_{j,i},s_{j,i})\in M_i\times [-2\epsilon_i^{-1},0]$ such that $(y_{j,i},s_{j,i})$ is  $(\epsilon_i,1)$-static and $(\epsilon_i,1)$-selfsimilar with $\mathcal{N}_{(y_i,s_i)}(1)\geq -Y$ and $\mathcal{N}_{(y_i,s_i)}(\epsilon_i)\leq \operatorname{log}(\frac{1}{2})+\epsilon_i$,  $(y_{2,i},s_{2,i})\in P^\ast(y_{1,i},0;\Lambda)\cap  P^\ast(y_{1,i},0;\gamma^3)^c$, but $|s_{2,i}|\leq c_i\gamma^6$. By \cite[Proposition 10.8, Proposition 12.1]{Bam3}, after relabeling our sequence, we may assume that $(y_{j,i},s_{j,i})$ is strongly $(1,\epsilon_i,1)$-split. In particular, by Theorem \ref{thm:bamconvergence} there is a correspondence $\mathfrak{C}$ such that  $$(M^4_i,(g_{i,t})_{t\in [-\epsilon_i^{-1},0]},(\nu_{y_{1,i},0;t})_{t\in [-\epsilon_i^{-1},0]})\xrightarrow[i\to\infty]{\mathbb{F},\mathfrak{C}} (\mathcal{X},(\nu_{y_{1,\infty};s})_{s\leq 0})$$ where $\mathcal{X}$ is a static flow modeled on a flat cone of the form $ C(\mathbb{S}^3/\Gamma)$, $|\Gamma|\leq e^Y$, with vertex $y_{1,\infty}$ such that $y_{1,i}\xrightarrow[i\to\infty]{\mathfrak{C}}y_{1,\infty}$. Moreover, we have that $C(\mathbb{S}^3/\Gamma)\cong X^\prime\times \mathbb{R}$ for some singular space $X^\prime$, and so $\Gamma$ must be trivial. This contradicts the fact that $\mathcal{N}_{(y_i,s_i)}(\epsilon_i)\leq \operatorname{log}(\frac{1}{2})+\epsilon_i$.

 	
 \end{proof}

The following lemma, which is the main result of this section, states that the center set of an $(m,\epsilon,\gamma,\Lambda)$-(quasi) annular region has bounded $2$-content, and that a small tubular neighborhood of the center set of an $(m,\epsilon,\gamma,\Lambda)$-\textit{annular} region contains all of the nearby 	almost singular points.

\begin{Lemma}[(Quasi) Annular Structure Lemma]\label{lem:annular}
    Let $\Lambda, Y<\infty$. Then for $\gamma\leq \overline{\gamma}(Y,\Lambda)$, and $\epsilon\leq \overline{\epsilon}(Y,\gamma,\Lambda)$, the following statement holds. Let $(M^4,(g_t)_{t\in [-2\epsilon^{-1},10]} )$ be a closed Ricci flow, $x_0\in M$ with $\mathcal{N}_{(x_0,0)}(1)\geq -Y$, and $\mathcal{A}=P^\ast(x_0,0;2\Lambda,2,2)\backslash \bigcup_{(x,t)\in\mathcal{C}} P^\ast(x,t;r_{(x,t)})$ be an $(m,\epsilon,\gamma,\Lambda)$-quasi annular region for some $m\in \mathbb{N}\backslash \{1\}$. Then we have the content estimate
    \begin{align}
    \sum_{(x,t)\in\mathcal{C}}r_{(x,t)}^2\leq C(Y,\gamma), \label{iden:annular_struc_1}
    \end{align}
     and if $\mathcal{A}$ is also an $(m,\epsilon,\gamma,\Lambda)$-annular region, we have in addition that for all $\sigma>0$, 
    \begin{align}
     \{(x,t)\in M \times [-1,1]\,\vert\, r_{\operatorname{Rm}}(x,t)<\sigma \}\cap \mathcal{A}\subseteq T_{D(Y,\gamma,\Lambda)\cdot\sigma}(\mathcal{C}), \label{iden:annular_struc_2}
    \end{align}
      where $T_{D\cdot \sigma}(\mathcal{C})$ is the tubular neighborhood of radius $D\sigma$ around $\mathcal{C}$ with respect to the natural topology.
\end{Lemma}

\begin{Remark} \label{rem:annular_struc_1}
	We will actually show the slightly stronger fact that for all $t\in [-1,1]$ $$ \operatorname{diam}_{P^\ast} \{(x,s)\in M \times [-1,1]\,\vert\, r_{\operatorname{Rm}}(x,s)<\sigma \}\cap \mathcal{A}\cap (M\times\{t\})\leq  D(Y,\gamma,\Lambda)\cdot \sigma.$$ 
\end{Remark}

\begin{proof}
    By adjusting $Y<\infty$ if necessary, we may assume that $\inf_{(x,t)\in P^\ast(x_0,0;2\Lambda)}\mathcal{N}_{(x,t)}(1)\geq -Y$.  For $(x,t), (x^\prime,t^\prime)\in \mathcal{C}$, set $\hat{r}((x,t),(x^\prime,t^\prime))=\max\{ r_{(x,t)}, r_{(x^\prime, t^\prime)}\}$. For $(x,t), (x^\prime,t^\prime)\in P^\ast(x_0,0;2\Lambda)$, let $\kappa((x,t),(x^\prime,t^\prime))>0$ be minimal such that $(x^\prime ,t ^\prime)\in P^{\ast}(x,t;\kappa)$.  We begin with the following claim.
    \begin{Claim}\label{claim:annular2}
For $\gamma\leq\overline{\gamma}(Y,\Lambda)$ and $\epsilon\leq \overline{\epsilon}(Y,\gamma,\Lambda)$, we have $\kappa\leq \hat{r}$ for every $(x,t),(x^\prime,t^\prime)\in \mathcal{C}$ with $|t-t^\prime|\leq \gamma^6\kappa^2$.
\end{Claim}
\begin{proof}
    Suppose the claim was false, so that there is a sequence $\epsilon_i\searrow 0$ together with a sequence of closed Ricci flows $(M^4_i,(g_{i,t})_{t\in [-2\epsilon_i^{-1},10]})$ and $(m_i,\epsilon_i,\gamma,\Lambda)$-quasi annular regions $\mathcal{A}^i\coloneqq P^\ast(x_{i,0},0;\Lambda,2,2)\backslash \bigcup_{(x,t)\in\mathcal{C}^i
    } P^\ast(x,t;r_{(x,t)}),$
    $\inf_{(x,t)\in P^\ast(x_{0,i},0;\Lambda)}\mathcal{N}_{(x,t)}(1)\geq -Y$, such that there exists $(y_i,s_i),(y_i^\prime,s_i^\prime)\in \mathcal{C}^i$ for which $$\kappa_i\coloneqq \kappa((y_i,s_i),(y_i^\prime,s_i^\prime))> \hat{r}_i\coloneqq \hat{r}_i((y_i,s_i),(y^\prime_i,s_i^\prime))$$ and $|s_i-s_i^\prime|\leq \gamma^6\kappa_i^2$. Assume without loss of generality that $s_i\geq s_i^\prime$. If $\limsup_{i\to\infty}|s_i-s_i^\prime|^{-1}=\infty$,  we are in the situation of the proof of Proposition \ref{prop:time_dif}, and so we can follow an identical blowup argument to get a contradiction. Therefore, we may assume that $\limsup_{i\to\infty}|s_i-s_i^\prime|^{-1}<\infty$.
   By Theorem \ref{thm:bamconvergence}, after rescaling by a factor of $\kappa_i$, there is a correspondence $\mathfrak{C}$ such that  we have the convergence
    $$ (M^4_i,(g_{i,t})_{t\in [-\epsilon_i^{-1},0]}, (\nu_{y_i,s_i;s})_{s\in [-\epsilon_i,0]}))\xrightarrow[i\to\infty]{\mathfrak{C},\mathbb{F}}(\mathcal{X}_\infty, \nu_{y_\infty;s})$$
    where $(\mathcal{X}_\infty, \nu_{y_\infty;s})$ is the static flow modeled on a flat cone of the form $C(\mathbb{S}^3/\Gamma)$, $|\Gamma|\leq e^Y$, with vertex $y_\infty$. Since $(y_i^\prime ,s_i ^\prime)\in P^\ast (y_i,s_i; \Lambda)$,
    we have using \cite[Lemma 4.1, Lemma 4.2]{fluck2025epsilonregularitytheoremnoncollapsedricci} that $(y_i^\prime,s_i^\prime)\xrightarrow[i\to\infty]{\mathfrak{C}} y^\prime_\infty$ for some $y^\prime_\infty\in \mathcal{S}_\infty$. Since $$|\mathfrak{t}(y_\infty^\prime)|= \lim_{i\to\infty} |s_i-s_i^\prime|\leq 2\gamma^6,$$
    we conclude that $y^\prime_\infty\in P^\ast(y_\infty; C(Y)\gamma)$. That is, $(y_i^\prime,s_i^\prime) \in P^\ast(y_i,s_i; 2C(Y)\gamma)$ for sufficiently large $i$. This contradicts the minimality of $\kappa$ if $\gamma \leq \overline{\gamma}(Y)$.
\end{proof}
Combining Claim \ref{claim:annular2} with Proposition \ref{prop:time_dif}, we get the following claim.
    \begin{Claim} \label{annular_claim1}
 For $\gamma\leq\overline{\gamma}(Y,\Lambda)$ and $\epsilon\leq \overline{\epsilon}(Y,\gamma,\Lambda)$, the following statement is true.  Suppose $(x,t), (x^\prime,t^\prime)\in \mathcal{C}$ satisfy $|t-t^\prime|< \gamma^7 \kappa^2$, then $(x,t)=(x^\prime,t^\prime)$.
\end{Claim}



Using Claim \ref{annular_claim1}, \eqref{iden:annular_struc_1} follows from an elementary real analysis argument. Now suppose that $\mathcal{A}$ is an $(m,\epsilon,\gamma,\Lambda)$-annular region.
By \eqref{iden:annular4} of Definition \ref{def:annular}, given $$(\tilde{x},\tilde{t})\in  \{(x,t)\in M\times [-1,1]\,\vert\, r_{\operatorname{Rm}}(x,t)<\sigma \}\cap \mathcal{A},$$ we can choose $(\tilde{y},\tilde{s})\in \mathcal{C}$ with $|\tilde{t}-\tilde{s}|\leq \gamma^5 r^2_{(\tilde{y},\tilde{s})}$. We have the following claim.
\begin{Claim}\label{claim:structure1}
    For $\gamma\leq \overline{\gamma}(Y,\Lambda)$, and $\epsilon\leq \overline{\epsilon}(Y,\gamma,\Lambda)$, we have $\kappa((\tilde{x},\tilde{t}),(\tilde{y},\tilde{s}))\leq D(Y,\gamma,\Lambda)\cdot\sigma$.
\end{Claim}
\begin{proof}
Suppose for the sake of contradiction that the claim was false. In particular, take some  $$(\tilde{x},\tilde{t})\in  \{(x,t)\in M\times [-1,1]\,\vert\, r_{\operatorname{Rm}}(x,t)<\sigma \}\cap \mathcal{A}$$ and $(\tilde{y},\tilde{s})\in \mathcal{C}$ with $|\tilde{t}-\tilde{s}|\leq \gamma^5 r^2_{(\tilde{y},\tilde{s})}$, and  suppose for the sake of contradiction that  $\kappa((\tilde{x},\tilde{t}),(\tilde{y},\tilde{s}))\geq D\cdot\sigma$. Note that by definition, $(\tilde{y},\tilde{s})$ is $(\epsilon,\kappa)$-static and $(\epsilon,\kappa)$-selfsimilar where $\kappa=\kappa((\tilde{x},\tilde{t}),(\tilde{y},\tilde{s}))$. Let $\theta,\alpha>0$ be constants to be determined. Take an $H_n$-center $(z,\tilde{s}-\kappa^2\theta)$ of $(\tilde{x},\tilde{t})$. 
Suppose that $r_{Rm}(z,\tilde{s}-\kappa^2\theta)\geq \alpha\kappa$. By \cite[Corollary 9.6]{Bam1}, we have that $$r_{Rm}(\tilde{x},\tilde{t})\geq c(\theta,\alpha)\kappa\geq Dc(\theta,\alpha)\sigma$$  if $\theta\leq \overline{\theta}(Y,\alpha)$, which is a contradiction if $D\geq \underline{D}(\theta,\alpha)$. Therefore, we may assume that $r_{Rm}(z,\tilde{s}-\kappa^2\theta)<\alpha\kappa$.  Choose $\theta\leq \overline{\theta}(Y)$ and  $\gamma\leq \overline{\gamma}(\theta)$ so that $\tilde{t}-\theta\kappa^2\in[\tilde{s}-2\theta\kappa^2, \tilde{s}-\frac{\theta}{2}\kappa^2]$ and  $(z,\tilde{s}-\kappa^2\theta)\in P^\ast(\tilde{y},\tilde{s};\frac{1}{2}\kappa)^c$. By smooth convergence on the regular set, and Perelman's forward pseudolocality Theorem \cite[Theorem 10.1]{perelman2002entropyformularicciflow}, we get that  $r_{Rm}(z,\tilde{s}-\kappa^2\theta)\geq c(Y)\kappa$ for $\epsilon\leq \overline{\epsilon}(Y,\theta,\gamma,\Lambda)$, from which we get a contradiction if $\alpha\leq \overline{\alpha}(Y)$.
\end{proof}
It is easy to see that Claim \ref{claim:structure1} implies \eqref{iden:annular_struc_2}, which completes the proof.
\end{proof}

The following proposition states that the time slices of an  $(m,\epsilon,\gamma,\Lambda)$-annular region are  close to the annulus of a flat cone in the $C^\infty$-topology.

\begin{Proposition}\label{prop:annular_time_slice}
     Let $\epsilon^\prime>0$ and $\Lambda,Y<\infty$. Then for $\gamma\leq \overline{\gamma}(Y,\epsilon^\prime,\Lambda)$ and $\epsilon\leq \overline{\epsilon}(Y,\gamma,\epsilon^\prime,\Lambda)$ the following statement holds. Let $(M^4,(g_t)_{t\in [-2\epsilon^{-1},10]} )$ be a closed Ricci flow, $x_0\in M$ with $\mathcal{N}_{x_0,0}(1)\geq -Y$, and $\mathcal{A}=P^\ast(x_0,0;2\Lambda,2,2)\backslash \bigcup_{(x,t)\in\mathcal{C}} P^\ast(x,t;r_{(x,t)})$ be an $(m,\epsilon,\gamma,\Lambda)$-annular region. Then there is some flat cone of the form $C(\mathbb{S}^3/\Gamma)$ where $|\Gamma|\leq e^{Y}$ such that the following holds. For each $t\in  [\inf\mathfrak{t}(P^\ast(x_0,t_0;1),\sup\mathfrak{t}(P^\ast(x_0,t_0;1)]$, there exists $r_{2,t}>r_{1,t}>0$, for which there is a diffeomorphism 
     $$\varphi_t : U_t\to A_{r_{1,t},r_{2,t}}$$
     where $ U_t\supseteq \mathcal{A}\cap (M\times \{t\})$ is open and  $A_{r_1,r_2}\subseteq C(\mathbb{S}^3/\Gamma)$ denotes the standard metric annulus, on which we have the estimate 
     \begin{align*}
         \| (\varphi_t^{-1})^\ast (g_{t})-g_{\operatorname{flat}}\|_{C^{\lfloor (\epsilon^\prime)^{-1}\rfloor }}\leq \epsilon^\prime
     \end{align*}
     where $g_{\operatorname{flat}}$ is the standard flat metric on $C(\mathbb{S}^3/\Gamma)$.
\end{Proposition}
\begin{proof}
By adjusting $Y$ if necessary, we may assume that $\inf_{(x,t)\in P^\ast(x_0,0;\Lambda)}\mathcal{N}_{(x,t)}(1)\geq -Y$.   Fix $t\in [-1,1]$ and let $\epsilon^{\prime\prime}>0$ be to be determined. By \eqref{iden:annular4} of Definition \ref{def:annular},   we can choose $(y_0,s_0)\in \mathcal{C}$ for which $|s_0-t|\leq 4\gamma^5 r_{(y_0,s_0)}^2$. 

By smooth convergence on the regular set, for all $c>0$, we can choose $\epsilon\leq \overline{\epsilon}(Y,\epsilon^{\prime\prime},\gamma,c,\Lambda)$, so that for all $r \in [r_{(y_0,s_0)},1]$, there exists a map 
    $$\varphi_{r}: A_{cr,2c^{-1}\Lambda r}\to P^\ast_{s_0-5\gamma^5r^2}(y_{0},s_{0};C(Y)c^{-1}\Lambda r)$$
 which is a diffeomorphism onto its image and for which we have    $$\| g_{\operatorname{flat}}-\varphi_{r}^\ast(g_{s_0-5\gamma^5r^2})\|_{C^{\lfloor (\epsilon^{\prime\prime})^{-1}\rfloor }}\leq \epsilon^{\prime\prime}.$$ By \cite[Corllary 9.6]{Bam1}, Perelman's forward pseudolocality Theorem \cite[Theorem 10.1]{perelman2002entropyformularicciflow}, and Shi's derivative estimates, if $\epsilon^{\prime\prime}\leq \overline{\epsilon^{\prime\prime}}(\epsilon^\prime,Y)$ and $\gamma \leq \overline{\gamma}(Y,\epsilon^\prime,\Lambda)$, we get a map
    $$\varphi_{t,r}: A_{a(Y)^{-1}\cdot cr,a(Y)\cdot 2c^{-1}\Lambda r}\to P^\ast_{t}(y_{0},s_{0};C(Y)c^{-1}\Lambda r)$$
     which is a diffeomorphism onto its image and for which we have    $\| g_{\operatorname{flat}}-\varphi_{t,r}^\ast(g_t)\|_{C^{\lfloor \epsilon^{-1}\rfloor }}\leq \epsilon^{\prime}.$ Since $$\mathcal{A} \subseteq P^\ast(y_0,s_0;4\Lambda)\backslash P^\ast(x_0,t_0;r_{(y_0,s_0)})$$ and $|s_0-t|\leq 2\gamma^5r_{(y_0,s_0)}^2$, we can choose $c\leq \overline{c}(Y)$ and $\gamma\leq \overline{\gamma}(Y,\epsilon^\prime ,\Lambda)$ so that $$\bigcup_{r\in [r_{(y_0,s_0)},1]}\varphi_{t,r}(A_{a(Y)\cdot cr,a(Y)\cdot 2c^{-1}\Lambda r})\supseteq \mathcal{A}\times (M\times \{t\}).$$  The proof now follows in the exact same way as \cite[Lemma 8.40]{regEinstein}.
\end{proof}
\section{Types of $P^\ast$-Parabolic Balls}\label{sec:ball}

In this section, we introduce the various types of $P^\ast$-parabolic balls that will play a role in the proof of Theorem \ref{thm:neck_decomp}.

We begin with two important results which motivate the precise formulation of Definition \ref{def:balls}. Our first lemma states that one cannot see arbitrarily small cylindrical necks for non-collapsed (real) $4$-dimensional \textit{K\"ahler} Ricci flows. Here the K\"ahler assumption is expected to be strictly necessary.

\begin{Lemma}[No Small Scale Cylindrical Necks] \label{lem:weaksplittingepsreg} For any $Y<\infty$, there exists $\epsilon_0=\epsilon_0(Y)>0$ such that the following holds. Suppose $(M^4,(g_t)_{t\in I})$ is a closed K\"ahler-Ricci flow, $r>0$, and $(x_0,t_0) \in M\times I$ is such that $[t_0-\epsilon_0^{-1}r^2,t_0+2r^2] \subseteq I$ and $\mathcal{N}_{x_0,t_0}(r^2)\geq -Y$. If $(x_0,t_0)$ is strongly $(1,\epsilon_0,r)$-split and $(\epsilon_0,r)$-selfsimilar, then $r_{Rm}(x_0,t_0)\geq \epsilon_0 r$. 
\end{Lemma}
\begin{proof} 
We may assume $t_0=0$ and $r=1$. Suppose by way of contradiction that the statement was false, then there is a sequence $\epsilon_i \searrow 0$ together with closed K\"ahler-Ricci flows $(M_i^4,(g_{i,t})_{t\in [-\epsilon_i^{-1},2]})$, points $x_{0,i} \in M_i$ such that $\mathcal{N}_{x_{0,i},0}(1)\geq -Y$, and $(x_{0,i},0)$ are strongly $(1,\epsilon_i,1)$-split and $(\epsilon_i,1)$-selfsimilar, but $r_{Rm}(x_{0,i},0)\leq \epsilon_i$. By Theorem \ref{thm:bamconvergence}, after passing to a subsequence, we have
$$(M_i,(g_{i,t})_{t\in [-\epsilon_i^{-1},0]},(\nu_{x_{i,0},0;t})_{t\in [-\epsilon_i^{-1},0]}) \xrightarrow[i\to \infty]{\mathbb{F}}(\mathcal{X},(\mu_t)_{t\in (-\infty,0)}),$$
on compact time intervals, where $(\mathcal{X},(\mu_t)_{t\in (-\infty,0)})$ is a metric soliton modeled on an orbifold shrinking gradient K\"ahler-Ricci soliton $(X,g,f)$ with at most isolated singularities, which splits as a metric product $X=X'\times \mathbb{R}$. Since the singularities of $X$ are isolated, it follows that $X$ is actually smooth, and because $X$ is K\"ahler, we have $X=X''\times \mathbb{C}$ as a metric product, where $X''$ is a shrinking K\"ahler-Ricci soliton of complex dimension 1.  Since $\lim_{i\to\infty }r_{Rm}(x_{0,i},0)=0$, we have $X''=\mathbb{CP}^1$ with the Fubini-Study metric. In particular, $X$ has a smooth holomorphic curve $\Sigma := \mathbb{CP}^1 \times \{0\}$ of self-intersection $0$. Let $\psi_i: U_i \to M_i$ be diffeomorphisms realizing the smooth convergence on the regular set, so that $\psi_i^{\ast}g_{i,-1} \to g$, $\psi_i^{\ast}J_i \to J$ in $C_{loc}^{\infty}(X)$, where $J_i,J$ are the complex structures of $M_i$ and $X$, respectively. By \cite{CCD}[Corollary 2.3], there is a (smooth) $\psi_i^{\ast}J_i$-holomorphic curve $\Sigma_i \subseteq M_i$ with self-intersection zero such that $\Sigma_i \subseteq B(\Sigma,1)$ and
$$\left| \int_{\Sigma_i} \psi_i^{\ast} \omega_{i,-1} -8\pi \right| = \left| \int_{\Sigma_i} \psi_i^{\ast} \omega_{i,-1} - \int_{\Sigma_i} \omega \right|<\pi$$
for sufficiently large $i\in \mathbb{N}$. Then $E_i :=\psi_i(\Sigma_i)$ is a $J_i$-holomorphic curve in $M_i$ with trivial self-intersection, and so we get that
$$\text{Vol}_{g_{i,t}}(E_i) = \int_{E_i} \omega_{i,t} = \int_{E_i} \omega_{i,-1} - 8\pi(t+1) \leq 9\pi -8\pi(t+1).$$
In particular, $\text{Vol}_{g_{i,1}}(E_i)<0$, contradicting the fact that $\omega_{i,1}$ is K\"ahler. 
\end{proof}

\begin{Corollary} \label{splittingepsreg} For any $Y<\infty$, there exists $\epsilon_0=\epsilon_0(Y)>0$ such that the following holds. Suppose $(M^4,(g_t)_{t\in I})$ is a closed K\"ahler-Ricci flow, $r>0$, and $(x_0,t_0) \in M\times I$ is such that $[t_0-2\epsilon_0^{-1}r^2,t_0+2r^2] \subseteq I$ and $\mathcal{N}_{x_0,t_0}(r^2)\geq -Y$. If $(x_0,t_0)$ is strongly $(1,\epsilon_0,r)$-split, then $r_{Rm}(x_0,t_0)\geq \epsilon_0 r$.     
\end{Corollary}
\begin{proof} We may assume $r=1$ and $t_0=0$. Let $\tilde{\epsilon}_0 = \tilde{\epsilon}_0(Y)>0$ be the constant from Lemma \ref{lem:weaksplittingepsreg}. By \cite[Proposition 7.1]{Bam3}, there is $s\geq r_0(Y,\tilde{\epsilon}_0) $ such that $(x_0,0)$ is $(\tilde{\epsilon}_0,s)$-selfsimilar. Because $(x_0,0)$ is strongly $(1,\epsilon_0,1)$-split, it is strongly $(1,\Psi(\epsilon_0|\tilde{\epsilon}_0,Y),s)$-split. If $\epsilon_0 \leq \overline{\epsilon}_0(\tilde{\epsilon}_0,Y)$, we may therefore apply Lemma \ref{lem:weaksplittingepsreg} to obtain $r_{Rm}(x_0,0) \geq s\tilde{\epsilon}_0 \geq r_0(Y,\tilde{\epsilon}_0)\epsilon_0$, from which the claim follows.
\end{proof}

Our next proposition states that if two almost self similar points are based at sufficiently different times, then they are also almost static.  This is a slight extension of \cite[Proposition 10.1]{Bam3}.

\begin{Proposition}\label{prop:max_sym}
	Let $\Lambda,Y<\infty$ and $\beta,\epsilon,\alpha >0$. For $\delta\leq \overline{\delta}(Y,\epsilon,\alpha,\beta,\Lambda)$, the following statement holds.  Let $(M^n,(g_t)_{t\in I} )$ be a closed Ricci flow, $x_0\in M$ be such that $\inf_{P^\ast(x_0,0;2\Lambda,2,2)}\mathcal{N}_{(x_0,0)}(1)\geq -Y$. If $(x_i,t_i)\in P^\ast(x_0,0;2\Lambda,2,2), i=1,2$, are  $(\delta,1)$-selfsimilar and $|t_1-t_2|\geq \alpha$, then
	\begin{enumerate}
		\item $(x_i,t_i)$ is $(\epsilon,s)$-static for all $s\in [\beta,1]$, \label{iden:ms1}
		\item there exists $m\in\mathbb{N}$ such that $|\mathcal{N}_{(x_i,t_i)}(\epsilon s^2)-\operatorname{log}(\frac{1}{m})|<\epsilon $ for all $s\in [\beta,\epsilon^{-1}]$. \label{iden:ms2}
	\end{enumerate}
\end{Proposition}
\begin{proof}
	The following claim is an immediate consequence of \cite[Proposition 7.1]{Bam3}.
	\begin{Claim}\label{claim:ms_1}
		For $\delta^\prime>0$ and $\delta \leq \overline{\delta}(Y,\beta,\delta^\prime)$, each $(x_i,t_i)$ is $(\delta^\prime,r)$-selfsimilar for each $r\in [\beta,1]$.
	\end{Claim}
	The following claim is a slight modification of \cite[Proposition 10.1]{Bam3}.
	\begin{Claim}\label{claim:ms_2}
		For $\delta\leq \overline{\delta}(Y,\epsilon,\alpha,\beta,\Lambda) $, each $(x_i,t_i)$ is $(\epsilon,r)$-static and $(\epsilon,r)$-selfsimilar for each $r\in [\beta,1]$.
	\end{Claim}
	\begin{proof}
		Fix $\delta^\prime>0$ to be determined.  By Claim \ref{claim:ms_1}, if $\delta \leq \overline{\delta}(Y,\beta,\delta^\prime)$, we can assume that $(x_i,t_i)$ is $(\delta^\prime,r)$-selfsimilar for each $i=1,2$ and $r\in [\beta,1]$.  Assume without loss of generality that $t_2\geq t_1$. Choose $\kappa(\beta,\alpha,\Lambda)>0$ such that $\kappa r^2\leq |t_2-t_1|\leq \kappa^{-1}r^2$ and $d_{W^1}^{g_{s}}(\nu_{(x_1,t_1);s},\nu_{(x_2,t_2);s})<\kappa^{-1}r^2$ for some $s\in [t_1-\kappa^{-1}r^2,t_1]$,  for each $r\in [\beta,1]$. It immediately follows from \cite[Proposition 10.1]{Bam3} that $(x_1,t_1)$ is $(\epsilon ,r)$-static for each $r\in [\beta,1]$ if $\delta^\prime \leq \overline{\delta^\prime}(\epsilon,Y,\kappa)$. Fix $r\in [\beta,1]$. As in the proof of \cite[Proposition 10.1]{Bam3}, for $\theta>0$, some $t^\ast\in [t_1-2\theta r^2,t_1-\theta r^2]$, and $t\in [t_1-\epsilon^{-1}r^2,t^\ast]$, we have that
		$$ \int _{M}\tau_1R\,d\nu_{(x_2,t_2);t}\leq C(Y,\kappa)\theta$$
		where $\tau_i=t_i-t$. It follows that for $t\in [t_1-2r^2,t_1-r^2]$ we have
		$$\int _{M}\tau_2R\,d\nu_{(x_2,t_2);t}\leq (\frac{\tau_2}{\tau_1})C(Y,\kappa)\theta\leq C(Y,\kappa)\theta,$$
		and so from \cite[Proposition 7.6]{kunikawa2025splittingquantitativestratificationsuper} (see also \cite[Discussion preceding Proposition 7.9]{Bam3}),
        we get that
        $$\int _{M}\tau_2R\,d\nu_{(x_2,t_2);t}\leq  C(Y,\kappa)\theta
		$$
        for all $t\in [t_2-\epsilon^{-1}r^2,t_2-\epsilon r^2]$ if $\delta^\prime\leq \overline{\delta^\prime}(\epsilon,Y,\kappa,\theta)$. In particular, for each $t\in [t_2-\epsilon^{-1},t_2-\epsilon]$ we have
        $$\int _{M}R\,d\nu_{(x_2,t_2);t}\leq \frac{\epsilon}{2}$$
		if $\theta\leq \overline{\theta}(Y,\kappa,\epsilon)$.  We get that 
		$$2r^2\int_{t_2-\epsilon^{-1}r^2}^{t_2-\epsilon r^2}\int_M |Rc|^2\,d\nu_{x_2,t_2;t}dt=\int_M R\,d\nu_{x_2,t_2;t}\arrowvert_{t_2-\epsilon^{-1}}^{t_2-\epsilon}\leq \epsilon$$
		for $\delta^\prime\leq \overline{\delta^\prime}(\epsilon,Y,\kappa,\theta)$. This shows that $(x_2,t_2)$ is $(\epsilon,r)$-static, which concludes the proof of the claim.
	\end{proof}
	Property \eqref{iden:ms1} follows immediately from Claim \ref{claim:ms_2}. Property \eqref{iden:ms2} follows easily from a compactness contradiction argument (c.f. Theorem \ref{thm:bamconvergence}).
\end{proof}


Based on Proposition \ref{prop:max_sym}, the following condition is natural.

\begin{Definition}\label{def:temp_ind}
    Let $\alpha,r>0$, and $(M^4,(g_t)_{t\in I} )$ be a closed Ricci flow. We say that  $Y\subseteq  M\times I$ is $(\alpha,r)$-temporally independent if there exists $(y_1,s_1), (y_2,s_2)\in Y$ with $|s_2-s_1|\geq \alpha \cdot r^2 $.
\end{Definition}

Our next condition guarantees that, for a given spacetime subset, the image under the projection onto the time factor is approximately dense in some interval  $J\subseteq \mathbb{R}$.

\begin{Definition}\label{def:temp_dens}
        Let $\eta,r>0$, and  $(M^4,(g_t)_{t\in I} )$ be a closed Ricci flow. We say that $Y\subseteq M\times I$ is $(\eta,r)$-temporally dense in $J\subseteq \mathbb{R}$ if $$T_{\eta\cdot r^2}(\mathfrak{t}(Y))\supseteq J,$$
        where $T_{\eta\cdot r^2}$ is the $(\eta\cdot r^2)$-tubular neighborhood in $\mathbb{R}$.
\end{Definition}

Next, we stratify $P^\ast$-parabolic balls  based on the behaviour of their entropy pinched points (c.f. \cite[Section 7.1]{L2CurvEin}).

\begin{Definition}\label{def:balls}
  Let $r,\eta, \delta,\alpha>0$, and $\Lambda,W<\infty$.  Let $(M^4,(g_t)_{t\in I} )$ be a closed Ricci flow, and $(x_0,t_0)\in M\times I$ be such that $\inf_{(x,t)\in P^\ast(x_0,t_0;\Lambda r)}\mathcal{N}_{x,t}(\delta^{-1}r^2)\geq -W-\delta$. Define 
    $$V_{\delta,W,r}(x_0,t_0)\coloneqq\{ (x,t)\in P^\ast(x_0,t_0;\Lambda r,2 r^2,2r^2)\,\vert\, |\mathcal{N}_{(x,t)}(\delta r^2)+W|\leq \delta \}.$$
    We say that $P^\ast(x_0,t_0;r)$ is 
    \begin{enumerate}
        \item a $(c)$-ball if $V_{\delta,W,r}$ is $(\eta,r)$-temporally dense in $[\inf\mathfrak{t}(P^\ast(x_0,t_0;r),\sup\mathfrak{t}(P^\ast(x_0,t_0;r)]$,
        \item an $(f)$-ball if $V_{\delta,W,r}$ is $(\alpha,r)$-temporally independent, but $V_{\delta,W,r}$ is not $(\eta,r)$-temporally dense in $[\inf\mathfrak{t}(P^\ast(x_0,t_0;r),\sup\mathfrak{t}(P^\ast(x_0,t_0;r)]$,
        \item a $(d)$-ball if $V_{\delta,W,r}$ is non-empty but not $(\alpha,r)$-temporally independent,
        \item an $(e)$-ball if $V_{\delta,W,r}$ is empty.
    \end{enumerate}
\end{Definition}

\begin{Remark}\label{rem:f_sub_int}
	If $P^\ast(x_0,t_0;r)$ is an $(f)$-ball, due to the lack of temporal density, we can choose a connected open subset of  $[\inf\mathfrak{t}(P^\ast(x_0,t_0;r)),\sup\mathfrak{t}(P^\ast(x_0,t_0;r))]$, which we denote by $I_{(x_0,t_0)}$, such that $|I_{(x_0,t_0)}|\geq \eta \cdot r^2$ and $P^\ast_{I_{(x_0,t_0)}}(x_0,t_0;r)\cap V_{\delta,W,r}(x_0,t_0)=\emptyset$. 
\end{Remark}

We finish this section by expanding on the meaning behind Definition \ref{def:balls} (see also \cite[Discussion on Page 75]{L2CurvEin}). By Proposition \ref{prop:max_sym} (see also \cite[Proposition 10.1]{Bam3}), a $(c)$-ball is one where we are able to construct $(m,\epsilon,\gamma)$-annular regions. By the same logic, an $(f)$-ball is one where we may only construct a quasi-static annular region. The presence of $(f)$-balls are technical obstructions in proving Theorem \ref{thm:diffeo} and Theorem \ref{thm:minkowski}, but we will see that we may decompose such balls into $(c),(d)$, and $(e)$-balls (c.f. Lemma \ref{lem:f_ball}), which, when combined with an additional induction argument, will allow us to remove all of the quasi-static regions from our covers (c.f. Lemma \ref{lem:ref_c} and Lemma \ref{lem:ref_d}). A $(d)$-ball is one in which we do not have any additional structure, but there is no strict entropy drop. These will be handled by showing that any such ball may be decomposed into $(e)$-balls, on which there is a strict entropy drop, and a collection of $(c)$-balls of very small content (c.f. Lemma \ref{lem:ref_d}). An $(e)$-ball is one on which there is a strict entropy drop. These will play a role in the final induction argument used to prove Theorem \ref{thm:neck_decomp}.



\section{The $0$-Dimensional Case}\label{sec:toy}

In this section, we prove the $0$-dimensional packing argument that will be needed in the proof of Theorem \ref{thm:neck_decomp}. One of the main new technical ingredients of this paper is Lemma \ref{lem:f_ball}. By Proposition \ref{prop:time_dif}, Proposition \ref{prop:max_sym}, and \cite[Proposition 10.8]{Bam3},
on an $(f)$-ball, the set $V_{\delta,W,1}$ of entropy pinched points will approximately form a temporal line. Away from a small tubular neighborhood of this line, one expects that any covering will contain only $(e)$-balls. Therefore, to produce the desired covering of our $(f)$-ball, it is essentially sufficient to construct an appropriate covering of an interval in $\mathbb{R}$. In this setting, our setup reduces to the following.

\begin{Definition}\label{def:frequency_function}
    A frequency function on an interval $I\subseteq \mathbb{R}$ is a continuous map $\mathcal{N}_t : I\times  \mathbb{R}_{\geq 0}\to \mathbb{R}_{\leq 0}$ such that:
    \begin{enumerate}
        \item $\mathcal{N}_t(0)=0$ for all $t\in I$,
        \item $\mathcal{N}_t=\mathcal{N}\vert_{\{t\}\times \mathbb{R}_{\geq 0 }}: \mathbb{R}_{\geq 0}\to \mathbb{R}_{\leq 0}$ is non-increasing for all $t\in I$.
    \end{enumerate}
\end{Definition}

The following definition will simplify notation substantially.

\begin{Definition}\label{def:int_ball}
    Given $t\in \mathbb{R}$, and $r,\gamma>0$, we set $J_{t,r,\gamma}=(t-\gamma\cdot r^2,t+\gamma\cdot r^2)$.
\end{Definition}

\begin{Definition}\label{def:int_types}
    Let $\mathcal{N}$ be a frequency function on the interval $(-2,2)$, $ \eta,\alpha ,\delta>0$,  and $t\in (-1,1)$ with $\inf_{(t-2r^2,t+2r^2)}\mathcal{N}_{t^\prime}(\delta^{-1}r^2)\geq -W-\delta$ for some $W<\infty$. Define 
$$V_{\delta,W,r}(t)\equiv \{ t^\prime\in J_{t,r,1}\, \arrowvert \, |\mathcal{N}_{t^\prime}(\delta  r^2)+W|<\delta\}.$$ We call the interval $J_{t,r,1}$ 
\begin{enumerate}
    \item a $(c)$-interval if for each $s\in J_{t,r,1}$ there exists $s^\prime \in V_{\delta,W,r}(t)$ such that $|s-s^{\prime}|<\eta r^2$,
    \item an $(f)$-interval if $V_{\delta,W,r}(t)$ if there exists  $s,s^\prime\in V_{\delta,m,r}(t)$  such that $|s-s^\prime|\geq \alpha r^2$, but there exists $s\in J_{t,r,1}$ such that for each $s^\prime \in V_{\delta,W,r}(t)$, we have $|s-s^{\prime}|\geq \eta r^2$,
    \item a $(d)$-interval if  $V_{\delta,W,r}(t)$ is non-empty but there does not exist $s,s^\prime\in V_{\delta,W,r}(t)$  such that $|s-s^\prime|\geq \alpha r^2$,
    \item an $(e)$-interval if $ V_{\delta,W,r}(t)$ is empty.
\end{enumerate}
As in Remark \ref{rem:f_sub_int}, if  $J_{t,r,1}$ is an $(f)$-interval, we let $I_t\subseteq J_{t,r,1}$ denote a connected open subset such that $|I_t|\geq \eta\cdot r^2$ and $V_{\delta,W,1}(t)\cap I_t=\emptyset $.
\end{Definition}

The following real analysis lemma states that an appropriately chosen collection of $(f)$-intervals must have uniformally bounded content. In Section \ref{sec:neck_decomp} we will argue that Lemma \ref{lem:f_ball} and Lemma \ref{lem:ref_c} may be reduced to the setting of this result.

\begin{Lemma}\label{lem:real_analysis}
	Let $\{T_s\}_{s\geq 0}$ be an increasing family of closed subsets of $(-4,4)$ and $\beta,\eta>0$. Then if $\{ J_{t_f,r_f;1}\}_{f\in F}$ is any countable collection of subintervals in $(-2,2)$ such that
	\begin{enumerate}
		\item for each $f\in F$, $r_f=\beta^j$ for some $j\in \mathbb{N}$,
		\item for each $f\in F$, $t_f\in T_{\beta^{-1}r_f}$,
        \item for each $f\in F$, there is a subinterval $I_f\subseteq J_{t_f,r_f;1}$ such that $|I_f|\geq \eta r_f^2$ and $I_f\cap T_r=\emptyset$ for any $r\leq r_f$,
		\item there is a uniform constant $M<\infty$ such that for each $k\in \mathbb{N}$, if $F^k\subseteq F$ is such that $r_f=\beta^k$, then at most $M$ elements of $F^k$ may intersect each other simultaneously,
	\end{enumerate}
	we must have $$\sum_{f\in F}r_f^2 \leq C(M,\eta).$$
\end{Lemma}
\begin{Remark}
    We will be most interested in the case where $T_s=\bigcup_{t\in (-2,2)}V_{\delta,W,s}(t)$ and the collection $\{ J_{t_f,r_f;1}\}_{f\in F}$ are $(f)$-intervals.
\end{Remark}

\begin{proof}
Define $G_j=(-2,2)\backslash T_{\beta^j}$, so that the sets $G_j$ are increasing and each component of $G_j$ is an open interval whose closure has end points in $T_{\beta^j}$. Given $t\in G_j$, let $l_j(t)$ denote the length of the connected component of $G_j$ containing $t$. Silghtly abusing notation, we let $l_j(t)=0$ for $t\in T_{\beta^j}$. Note that for any fixed $t\in (-2,2)$, $l_j(t)$ is increasing in $j$. Set $A_k=\bigcup_{f\in F^k}I_{f}$. We have 
\begin{align}\label{iden:real_anal_1}
	\sum_{f\in F}r_f^2 \leq \eta^{-1}\sum_{k\in \mathbb{N}}\sum_{f\in F^k}|I_f|\leq 4M\eta^{-1}\sum_{k\in \mathbb{N}}|A_k|.
\end{align}
Write $A_k$ in the form $A_k=A^\prime_k\cup A_k^{\prime\prime}$ where $A_k^\prime=\{t\in A_k\,\vert\, l_{k-1}(t)\geq \eta\beta^{2k} \}$. We will bound the measure of $A_k^\prime$ and $A_k^{\prime\prime}$ seperately. First we turn to $A_k^{\prime\prime}$. Suppose that there exist $t\in A^{\prime \prime }_j\cap A_k$ for some $j<k$. It follows that $$l_{k-1}(t)<\eta\beta^{2k}<\eta\beta^{2j}\leq l_{j}(t)\leq l_{k-1}(t) ,$$ 
which is a contradiction. Therefore, the sets $A_k^{\prime\prime}$ are pairwise disjoint, and 
\begin{align}\label{iden:real_anal_2}
	\sum_{k\in \mathbb{N}}|A_k^{\prime\prime}|\leq 4.
\end{align}
Next we turn to $A_k^{\prime}$. Consider some  $t\in A_k^\prime$ and suppose that $t\in I_f$ for some $f\in F^k$. Since $t_f\in V_{\delta ,W,\beta^{-1}r_f}$, we must have that the connected component of $G_{k-1}\cap A_k^\prime$ containing $t$ is of length at most $2 \beta^{2k}$. In particular, it follows that 
$$|A_k^\prime| \leq \int_{-2}^2 w_k\,dt $$
where $w_k\coloneqq \frac{2\beta^{2k}}{l_{k-1}(t)}\chi_{\{t\in (-2,2)\, \vert\, l_{k-1}(t)\geq \eta\beta^{2k}  \}}$. Fix $t\in (-2,2)$. Suppose that $\{k \in \mathbb{N}\, \vert \, l_{k-1}(t)\geq \eta\beta^{2k}\}$ is non-empty and let $k_0$ be the minimum value. It follows that $$\sum_{k\in \mathbb{N}}w_k(t)\leq \frac{1}{\eta}\sum_{k\geq k_0}\beta^{2(k-k_0)}\leq C(\eta).$$
In particular, it we get that
\begin{align}\label{iden:real_anal_3}
\sum_{k\in\mathbb{N}}|A_k^\prime|\leq \int_{-2}^2w_k(t)\,dt \leq C(\eta).
\end{align}
Combining \eqref{iden:real_anal_1}-\eqref{iden:real_anal_3} easiy gives the result.
	\end{proof}

We end this section by proving an $(f)$-interval covering lemma. This states that any $(f)$-interval may be decomposed into $(c), (d),$ and $(e)$ intervals, whilst retaining control on the $2$-content (c.f. Lemma \ref{lem:f_ball}).  The spacetime version, Lemma \ref{lem:f_ball}, will follow essentially the same argument.

\begin{Lemma}
  Let $\mathcal{N}$ be a frequency function on $(-4,4)$, $ \delta,\eta,\alpha>0$, and suppose that $\inf_{t\in (-2,2)}\mathcal{N}_t(\delta^{-1})\geq -W-\delta$ for some $\delta<W<\infty$. If $J_{0,1,1}$ is an $(f)$-interval, there exists a covering 
    $$J_{0,1,1}\subseteq \bigcup_{c\in C}J_{t_c,r_c,1}\cup \bigcup_{d\in D}J_{t_d,r_d,1}\cup \bigcup_{e\in E}J_{t_e,r_e,1}$$
    by $(c),(d),$ and $(e)$-intervals, for which we have the content estimate
    \begin{align*}
        \sum_{c\in C}r_c^2+\sum_{e\in E}r_e^2+\sum_{d\in D}r_d^2\leq C(\eta).
    \end{align*}
\end{Lemma}
\begin{proof}
  Let $\gamma,\beta>0$ be constants to be determined. Choose a maximal subsets $\{t_i\}_{i\in \tilde{N}^1}\subseteq V_{\delta ,W,1}(0)$ such that $\{J_{t_i,\beta,\gamma^6}\}_{i\in \tilde{N}^1}$ are pairwise disjoint, and a collection and $\{t_i\}_{i\in \hat{N}^1}\subseteq (-1,1)\backslash T_{10\gamma^3\beta}(V_{\delta,W,1}(0))$ such that $\{J_{t_i,\beta,\frac{\gamma^6}{40}}\}_{i\in \hat{N}_1}$  are pairwise disjoint,
    
 Set $r_i=\beta$ for $i\in \tilde{N}^1$ and $r_i=\frac{\gamma^3}{40}\beta$ for $i\in\hat{N}^1$. Write $N^1\equiv \hat{N}^1\cup \tilde{N}^1=C^1\cup D^1\cup E^1\cup F^1$ so that:
    \begin{enumerate}
        \item  for each $c\in C^1$, $J_{t_c,r_c,1}$ is a $(c)$-interval,
        \item  for each $d\in D^1$, $J_{t_d,r_d,1}$ is a $(d)$-interval,
        \item  for each $e\in E^1$, $J_{t_e,r_e,1}$ is a $(e)$-interval,
        \item  for each $f\in F^1$, $ J_{t_f,r_f,1}$ is an $(f)$-interval.
    \end{enumerate}
By construction, we must have $\hat{N}^1=E^1$, so that any $(c)$ or $(f)$-interval in our cover is centered at a point in $V_{\delta,W,1}(0)$.

 For each $f\in F^1$, $J_{t_f,r_f,1}$ satisfies the hypothesis of the lemma, so we may repeat the construction on these intervals. That is, we choose a maximal subset  $\{t_i\}_{i\in \tilde{N}^2}\subseteq \bigcup_{f\in F^1}V_{\delta,W,\beta}(t_f)$ such that $\{J_{t_i,\beta\cdot r_f,\gamma^6}\}_{i\in \tilde{N}^2_f}$ are pairwise disjoint, and a maximal subset  $$\{t_i\}_{i\in \hat{N}^2}\subseteq \bigcup_{f\in F^1}J_{t_f,r_f;1}\backslash T_{\frac{\gamma^3}{40}\beta^2}\left(\bigcup_{f\in F^1}V_{\delta,W,\beta}(t_f)\right)$$ such that 
 $\{J_{t_i,\beta\cdot r_f,\frac{1}{40}\gamma^6}\}_{i\in \hat{N}^2}$ are pairwise disjoint.
 
 
 Set $r_i=\beta \cdot r_f$ for $i\in \tilde{N}^2$ and $r_i=\frac{\gamma^3}{40}\beta\cdot r_f$ for $i\in \hat{N}^2$, and write $N^2=C^2\cup D^2\cup E^2\cup F^2$ so that:
    \begin{enumerate}
        \item  for each $c\in C^2$, $(t_c-r_c^2,t_c+r_c^2)$ is a $(c)$-interval,
        \item  for each $d\in D^2$, $(t_d-r_d^2,t_d+r_d^2)$ is a $(d)$-interval,
        \item  for each $e\in E^2$, $(t_e-r_e^2,t_e+r_e^2)$ is a $(e)$-interval,
        \item  for each $f^\prime \in F^2$, $(t_f-r_{f^\prime}^2,t_f+r_{f^\prime}^2)$ is a $(f)$-interval.
    \end{enumerate}
    Set
    \begin{align*}
        \begin{cases}
            C^2& = C^1\cup C^2 \\  D^2& = D^1\cup D^2 \\ E^2& =E^1 \cup  E^2 
        \end{cases},
    \end{align*}
    so that we now have a cover
    \begin{align*}
        J_{0,1,1}&\subseteq \bigcup_{c\in C^2}J_{t_c,r_c,1}\cup \bigcup_{d\in D^2}J_{t_d,r_d,1}  \cup \bigcup_{c\in E^2}J_{t_e,r_e,1}\cup\bigcup_{f\in F^2}J_{t_f,r_f,1}
          \end{align*}
   where $r_f=\beta^2$ for each $f\in F^2$. Note that by construction, we have $J_{t_f,r_f,\gamma^6}\cap J_{t_{f^\prime},r_{f^\prime},\gamma^6}=\emptyset$ if $r_f=r_{f^\prime}$.

     Repeating this procedure $k$-times we get a sequence of coverings
       \begin{align*}
        J_{0,1,1}&\subseteq \bigcup_{c\in C^l}J_{t_c,r_c,1} \cup \bigcup_{d\in D^l}J_{t_d,r_d,1}   \cup \bigcup_{c\in E^l}J_{t_e,r_e,1}\cup\bigcup_{f\in F^l}J_{t_f,r_f,1}
          \end{align*}
         for $l\leq k$ where $r_f=\beta^l$ for each $f\in F^l$,
         \begin{enumerate}
         \item for each $ f,f^\prime \in \bigcup_{l=1}^k F^l$ we have $J_{t_f,r_f,\gamma^6}\cap J_{t_{f^\prime},r_{f^\prime},\gamma^6}=\emptyset$ if $r_f=r_{f^\prime}$,
         \item for each $f\in F^k$ we have $t_f\in V_{\delta,W,\beta^{-1}r_{f}}$.
         \end{enumerate}
     By Lemma \ref{lem:real_analysis}, we get that 
         $$\sum_{l=0}^k\sum_{f\in F^k}r_f^2\leq C(\eta).$$
         In particular,
                   \begin{align*}
              \sum_{l=1}^{k}\sum_{c\in C^l}r_c^2+ \sum_{l=1}^{k}\sum_{d\in D^l}r_d^2+  \sum_{l=1}^{k}\sum_{e\in E^l}r_e^2\leq C_0^\prime\sum_{l=1}^{k-1}\sum_{f\in F^l}r_f^2\leq C(\eta).
          \end{align*}
  By continuity of the frequency function, there is some $k_0\in \mathbb{N}$ such that the $k_0$-th element of the above sequence  contains no $(f)$-intervals. It is easy to see that this covering satisfies the desired properties, which completes the proof.
\end{proof}

\section{Construction of the Neck Decomposition}\label{sec:neck_decomp}

In this section, we prove Theorem \ref{thm:neck_decomp}. First, we do an initial covering of each ball type from Definition \ref{def:balls}.

\begin{Lemma}[Coarse $(c)$-ball covering] \label{lem:coarse_c}
      Let $Y,\Lambda <\infty$ and $\beta,\sigma, \alpha, \gamma,\epsilon >0$. For $\eta\leq \overline{\eta}(Y,\gamma,\beta)$ and $\delta \leq \overline{\delta}(Y,\sigma,\epsilon,\gamma,\eta,\alpha,\Lambda,\beta)$,  the following statement is true.  Let $(M^4,(g_t)_{t\in [-2\delta^{-1},10]} )$ be a closed K\"ahler Ricci flow and $x_0\in M$ be such that $ \inf_{P^\ast(x_0,0;2\Lambda)}\mathcal{N}_{(x,t)}(\delta^{-1})\geq -W-\delta$ for some $\sigma \leq W\leq Y-\delta$. Let $\bigcup_{f\in \tilde{F}}P^\ast(x_f,t_f;r_f)\subseteq P^\ast(x_0,0;2\Lambda,2,2)$ be a family of $(f)$-balls such that
      \begin{enumerate}
          \item for every $f\in \tilde{F}$, $r_f=\beta^k$ for some $k\in \mathbb{N}$,
          \item for every $f\in \tilde{F}$, $(x_f,t_f)\in V_{\delta ,W,\beta^{-1}r_f}(x,t)$ for some $(x,t)\in P^\ast(x_0,0;2\Lambda)$,
          \item  for every $f,f^\prime \in \tilde{F}$, we have $P^\ast(x_f,t_f;\gamma^3r_f)\cap P^\ast(x_{f^\prime},t_{f^\prime};\gamma^3r_{f^\prime})=\emptyset $ if $r_f=r_{f^\prime}$.
      \end{enumerate}
  If $P^\ast(x_0,0;1)$ is a $(c)$-ball, then for some $m\leq e^{Y}$  there exists a covering
    $$ P^\ast(x_0,0;1)\subseteq \mathcal{A}\cup \bigcup_{d\in D}P^\ast(x_d,t_d;r_d)\cup \bigcup_{e\in E}P^\ast(x_e,t_e;r_e)\cup \bigcup_{f\in F}P^\ast(x_f,t_f;r_f)$$
    where 
   \begin{enumerate}
       \item $\mathcal{A}=P^\ast(x_0,0;2\Lambda,2,2)\backslash\bigcup_{i\in\mathcal{C}}P^\ast(x_i,t_i;r_i)$, $\mathcal{C}\coloneqq D\cup E\cup F$, is an $(m,\epsilon,\gamma,\Lambda)$-annular region,
       \item for every $d\in D$, $P^\ast(x_d,t_d;r_d)$ is a $(d)$-ball,
       \item for every $e\in E$, $P^\ast(x_e,t_e;r_e)$ is an $(e)$-ball, 
       \item for every $f\in F$, $P^\ast(x_f,t_f;r_f)$ is an $(f)$-ball.

   \end{enumerate}
    Moreover, all properties of the collection $\tilde{F}$ continue to hold for the collection $F\cup \tilde{F}$.
\end{Lemma}

\begin{Remark}
    The collection $\bigcup_{f\in \tilde{F}}P^\ast(x_f,t_f;r_f)$ is a purely technical feature that is required in order to repeatedly apply Lemma \ref{lem:coarse_c} in the proof of Lemma \ref{lem:ref_c}. When first reading, it may be instructive to go through the proof with $\tilde{F}=\emptyset$. 
\end{Remark}

\begin{Remark}
	Technically, we only include  $\mathcal{A}\cap (M\times [\inf \mathfrak{t}(P^\ast(x_0,0;1)),\sup \mathfrak{t}(P^\ast(x_0,0;1))])$ in our cover, as this is where we have control over the geometry of the time slices of $\mathcal{A}$ (c.f. Proposition \ref{prop:annular_time_slice}). To simplify our expressions, we ignore this fact in our notation.
	\end{Remark}



\begin{proof}
   By Proposition \ref{prop:max_sym} we may assume that $W=|\operatorname{log}(\frac{1}{m})|$ for some $m\in \mathbb{N}$. Choose a maximal subset $$\{ (x_i,t_i)\}_{i\in N^1}\subseteq \left(P^\ast(x_0,0;\Lambda,1,1)\cap V_{\delta,W,1}(x_0,0)\right)\backslash \bigcup_{f\in \tilde{F}, r_f=\beta}P^\ast(x_f,t_f;4\gamma^3\beta)$$ such that $P^\ast(x_i,t_i;\gamma^3\cdot \beta)$ are pairwise disjoint. Write $\hat{N}^1=C^1\sqcup D^1\sqcup E^1\sqcup \hat{F}^1$ and $r_i\coloneqq \beta$ where 
    \begin{enumerate}
        \item for each $c\in C^1$, $P^\ast(x_c,t_c;r_c)$ are $(c)$-balls,
        \item for each $d\in D^1$, $P^\ast(x_d,t_d;r_d)$ are $(d)$-balls,
        \item for each $e\in E^1$, $P^\ast(x_e,t_e;r_e)$ are $(e)$-balls,
       \item  for each $f\in \hat{F}^1$, $P^\ast(x_f,t_f;r_f)$ are $(f)$-balls.
    \end{enumerate}
  Let $\tilde{F}^1\subseteq \tilde{F}$ be such that for each $f\in \tilde{F}^1$ we have $r_f=\beta$ and 
  $$P^\ast(x_f,t_f;4\gamma^3r_f)\cap V_{\delta,W,1}(x_0,0)\cap P^\ast(x_0,0;2\Lambda,2,2)\neq\emptyset,$$
  and write $N^1=\hat{N}^1\cup \tilde{F}^1$, $F^1=\hat{F}^1\cup \tilde{F}^1$.
 Since $V_{\delta,m,1}(x_0,0)$ is $(\eta,1)$-temporally dense in $[\inf\mathfrak{t}(P^\ast(x_0,0;1),\sup\mathfrak{t}(P^\ast(x_0,0;1)]$, we have 
    \begin{align*}
          \bigcup_{i\in N^1}J_{t_i,\beta,\gamma^5}&\supseteq \bigcup_{i\in N^1}J_{t_i,\beta, 100\gamma^6+\eta\beta^{-2}}\\&\supseteq T_{\eta}(\mathfrak{t}\bigcup_{i\in N^1}P^\ast(x_i,t_i; 10\gamma^3\beta))\\&\supseteq T_{\eta}(\mathfrak{t}(V_{\delta,m,1}(x_0,0))\\&\supseteq [\inf\mathfrak{t}(P^\ast(x_0,0;1)),\sup\mathfrak{t}(P^\ast(x_0,0;1))] 
    \end{align*}
  for $\gamma\leq \gamma_0$ and $\eta\leq \overline{\eta}(\gamma,\beta)$. Using this fact together with Proposition \ref{prop:max_sym} we have that $$\mathcal{A}^1=P^\ast(x_0,0;2\Lambda,2,2)\backslash\bigcup_{i\in N^1}P^\ast(x_i,t_i;r_i) $$ is an $(m,\epsilon,\gamma,\Lambda)$-annular region for  $\eta\leq \overline{\eta}(\gamma,\beta)$ and $\delta\leq \overline{\delta}(Y,\epsilon, \beta,\gamma,\eta,\alpha,\sigma)$.  
  
  
  

Next, choose a maximal subset 
$$\{(x_i,t_i)\}_{i\in \tilde{N}^{{2}}}\subseteq \bigcup_{c\in C^1}V_{\delta,m,r_c}(x_c,t_c)\backslash \left(\bigcup_{i\in N^1\backslash C^1}P^\ast(x_i,t_i;\gamma^3r_{i}+4\gamma^3\beta^{2})\cup \bigcup_{f\in \tilde{F}, r_f=\beta^2}P^\ast(x_f,t_f;4\gamma^3\beta^2)\right)$$
such that $\{P^\ast(x_i,t_i;\gamma^3\beta^{2})\}_{i\in \tilde{N}^{2}}$ are pairwise disjoint.    
Write $\tilde{N}^{2}=C^{2}\sqcup D^{2}\sqcup E^{2}\sqcup \hat{F}^{2}$ and $r_i\coloneqq \beta^{2}$ where 
\begin{enumerate}
	\item for each $c\in C^{2}$, $P^\ast(x_c,t_c;r_c)$ is a $(c)$-ball, 
	\item for each $d\in D^{2}$, $P^\ast(x_d,t_d;r_d)$ is a $(d)$-ball,
	\item for each $e\in E^{2}$, $P^\ast(x_e,t_e;r_e)$ is an  $(e)$-ball,
	\item  for each $f\in \hat{F}^{2}$, $P^\ast(x_f,t_f;r_f)$ is an $(f)$-ball.
\end{enumerate}
Let $\tilde{F}^2\subseteq \tilde{F}$ be such that for each $f\in \tilde{F}^2$ we have $r_f=\beta^2$ and
\begin{enumerate}
    \item $P^\ast(x_f,t_f;4\gamma^3r_f)\cap \bigcup_{c\in C^1}V_{\delta,W,\beta}(x_c,t_c)\neq\emptyset$,
    \item for each $i\in \tilde{N}^2\cup N^1\backslash C^1$, we have $P^\ast(x_f,t_f;4\gamma^3r_f)\cap P^\ast(x_i,t_i;4\gamma^3r_i)=\emptyset$.
\end{enumerate}
Write 
\begin{align*}
\begin{cases}
     N^{2}&=\tilde{N}^{2}\cup N^1\backslash C^1\cup \tilde{F}^2\\
     D^{2}&=D^{2}\cup D^1 \\ E^{2}&=E^{2}\cup E^1 \\ F^{2}&=\hat{F}^{2}\cup F^1\cup \tilde{F}^2
\end{cases}.
\end{align*}
 By  Proposition \ref{prop:max_sym}, properties \eqref{iden:annular1}-\eqref{iden:annular3} of Definition \ref{def:annular} are satisfied by $$\mathcal{A}^{2}= P^\ast(x_0,0;2\Lambda,2,2)\backslash \bigcup_{i\in N^{2}}P^\ast(x_i,t_i;r_i).$$ Moreover, 
\begin{align*}
	\bigcup_{i\in \tilde{N}^{2}}J_{t_i,\beta^{2},\gamma^5}&\supseteq \bigcup_{i\in \tilde{N}^{2}}J_{t_i,\beta^{2},100\gamma^6+\eta\beta^{-2}}\\  &\supseteq T_{\eta\beta^{4}} (\mathfrak{t} \bigcup_{i\in \tilde{N}^{2}}P^\ast(x_i,t_i;10\gamma^3\beta^{2}))\\&\supseteq  \bigcup_{c\in C^1} (t_c-r_{c}^2,t_c+r_{c}^2)\backslash \bigcup_{i\in N^1\backslash C^1\cup \tilde{F}^2}J_{t_i, r_i,\gamma^6+(\eta+\gamma^3\beta)\beta^{4}r_i^{-2}},
\end{align*}
and so we have 
\begin{align*}
	\bigcup_{i\in N^{2}}J_{t_i,r_i,\gamma^5}&\supseteq \bigcup_{i\in N^1\backslash C^1\cup \tilde{F}^2}J_{t_i,r_i,\gamma^5}\cup  \bigcup_{i\in \tilde{N}^{2}}J_{t_i,\beta^{2},\gamma^5}\\ & \supseteq \bigcup_{i\in N^1\backslash C^1\cup \tilde{F}^2}J_{t_i,r_i,\gamma^5}\cup    \left(\bigcup_{c\in C^1}(t_c-r_{c}^2,t_c+r_{c}^2)\backslash \bigcup_{i\in N^1\backslash C^1\cup \tilde{F}^2}J_{t_i, r_i,\gamma^6+(\eta+\gamma^3\beta)\beta^{4}r_i^{-2}}\right)\\& \supseteq \bigcup_{i\in N^1} J_{t_i,r_i,\gamma^5}\\ & \supseteq [\inf\mathfrak{t}(P^\ast(x_0,0;1)),\sup\mathfrak{t}(P^\ast(x_0,0;1))]
\end{align*}
for $\gamma \leq \gamma_0$ and $\eta\leq \overline{\eta}(\beta,\gamma)$. In particular, property \eqref{iden:annular4} of Definition \ref{def:annular} is satisfied, and so $\mathcal{A}^{2}$ is an $(m,\epsilon,\gamma,\Lambda)$-annular region. 
Iterating this procedure, we arrive at the following claim.

    \begin{Claim}\label{claim:c_ball_1}
        There exists a sequence of coverings
        $$ P^\ast(x_0,t_0;1)\subseteq \mathcal{A}^k\cup \bigcup_{c\in C^k}P^\ast(x_c,t_c;r_c)\cup  \bigcup_{d\in D^k}P^\ast(x_d,t_d;r_d) \cup  \bigcup_{e\in E^k}P^\ast(x_e,t_e;r_e) \cup  \bigcup_{f\in F^k}P^\ast(x_f,t_f;r_f),$$
        where $r_c=\beta^k$ for all $c\in C^k$ and if  $N^k=C^k\sqcup D^k \sqcup E^k\sqcup F^k$, then $$\mathcal{A}^k= P^\ast(x_0,0;2\Lambda,2,2)\backslash \bigcup_{i\in N^k}P^\ast(x_i,t_i;r_i) $$
        is an $(m,\epsilon, \gamma,\Lambda)$-annular region, and
         \begin{enumerate}
        \item for each $c\in C^k$, $P^\ast(x_c,t_c;r_c)$ is a $(c)$-ball,  \label{iden:c_ball_c1}
        \item for each $d\in D^k$, $P^\ast(x_d,t_d;r_d)$ is a  $(d)$-ball, \label{iden:c_ball_c2}
        \item for each $e\in E^k$, $P^\ast(x_e,t_e;r_e)$ is an $(e)$-ball, \label{iden:c_ball_c3}
       \item  for each $f\in F^k$, $P^\ast(x_f,t_f;r_f)$ is an $(f)$-ball.
    \end{enumerate}
    Moreover, the properties outlined in Lemma \ref{lem:coarse_c} for the collection $\tilde{F}$ are retained for the collection $\tilde{F}\cup F^k$.
    \end{Claim}
    Since $(M^4,(g_t)_{t\in [-2\delta^{-1},10]})$ is smooth and compact, there is some $k_0\in \mathbb{N}$ for which the $k_0$-th element of the sequence constructed in Claim \ref{claim:c_ball_1} contains no $(c)$-balls. Clearly this covering satisfies all the requirements of the lemma, which completes the proof.
\end{proof}

\begin{Lemma}[Coarse $(d)$-ball covering] \label{lem:coarse_d}
    Let $\eta,\sigma>0$ and  $Y,\Lambda<\infty$. Then for  $\alpha\leq \overline{\alpha}(Y,\Lambda)$ and  $\delta \leq \overline{\delta}(Y,\alpha,\sigma,\eta,\Lambda)$, the following statement holds. Let $(M^4,(g_t)_{t\in [-2\delta^{-1},10]} )$ be a closed K\"ahler Ricci flow, and $x_0\in M$ be such that $\inf_{P^\ast(x_0,0;2\Lambda)}\mathcal{N}_{(x_0,0)}(\delta^{-1})\geq -W-\delta$ for some $\sigma \leq W\leq Y-\delta$. If $P^\ast(x_0,0;1)$ is a $(d)$-ball, then there exists a covering
    $$ P^ \ast(x_0,t_0;1)\subseteq \bigcup_{c\in C}P^\ast(x_c,t_c;r_c)\cup \bigcup_{e\in E}P^\ast(x_e,t_e;r_e)\cup\bigcup_{f\in F}P^\ast(x_f,t_f;r_f)$$
    where 
    \begin{enumerate}
          \item for every $c\in C$, $P^\ast(x_c,t_c;r_c)$ is a $(c)$-ball,
       \item for every $e\in E$, $P^\ast(x_e,t_e;r_e)$ is an $(e)$-ball, \label{iden:d_ball_2}
       \item for every $f\in F$, $P^\ast(x_f,t_f;r_f)$ is an $(f)$-ball,
    \end{enumerate}
    and we have the content estimates
    \begin{align*}
       \sum_{f\in F}r_f^2+ \sum_{c\in C}r_c^2  &\leq C(Y,\Lambda)\cdot\alpha^{\frac{1}{4}} \\ \sum_{e\in E}r_e^2&\leq C(Y,\alpha,\Lambda) .
    \end{align*}
\end{Lemma}


\begin{proof}
Let $\kappa>0$ be a constant to be determined. By Lemma \ref{lem:weaksplittingepsreg}, for $\kappa \leq \overline{\kappa}(Y,\sigma)$,
any point $(x,t)\in P^{\ast}(x_{0},0;4)$ which is strongly
$(1,\kappa,r)$-split, $r\in(0,1]$, satisfies $r_{Rm}(x,t)\geq \sigma r$
and $\mathcal{N}_{(x,t)}(\kappa \cdot r^2)\geq -\frac{1}{4}\sigma$. Let $\lambda_{0}=\lambda_{0}(\kappa)>0$
be as in \cite[Proposition 10.8]{Bam3}, and choose $\alpha\leq\overline{\alpha}(Y)$
so that $\alpha\leq \beta D_{0}^{-4}\lambda_{0}^{4}$, where $\beta>0$ and  $D_{0}=D_{0}(Y)\geq2$ are also from \cite[Proposition 10.8]{Bam3}. Also set $\lambda:=(\frac{\alpha}{\beta})^{\frac{1}{4}}$,
so that $\lambda\leq\lambda_{0}$. 

Let $\{(x_{i},t_{i})\}_{i\in N^1}\subseteq P^{\ast}(x_{0},0;1)$ be maximal such that $\{P^{\ast}(x_{i},t_{i};\frac{\sqrt{\beta}}{10}\lambda)\}_{i\in N^1}$
are pairwise disjoint.  Write $N^1=C^1\sqcup D^1\sqcup E^1\sqcup F^1$, and $r_i=\frac{\sqrt{\beta}}{2}\lambda$, so that 
\begin{enumerate}
	\item for each $c\in C^{1}$, $P^\ast(x_c,t_c;r_c)$ is a $(c)$-ball, 
	\item for each $d\in D^{1}$, $P^\ast(x_d,t_d;r_d)$ is a $(d)$-ball,
	\item for each $e\in E^{1}$, $P^\ast(x_e,t_e;r_e)$ is an  $(e)$-ball, 
	\item  for each $f\in F^{1}$, $P^\ast(x_f,t_f;r_f)$ is an $(f)$-ball.
\end{enumerate}
 By \cite[Theorem 9.11]{Bam1}, we have $|E^{1}|\leq C(Y)\lambda^{-6}$, so that
$$\sum_{e\in E^{1}}r_{e}^{2}\leq C_{b}(Y)\alpha^{-2}.$$
\begin{Claim} \label{claim:d_ball}
  For   $\alpha\leq \overline{\alpha}(Y)$ and  $\delta \leq \overline{\delta}(Y,\alpha,\sigma,\eta,\Lambda)$, there exists $C(Y)<\infty$ such that $|C^1|+|D^1|+|F^1|\leq C(Y,\Lambda)$.
\end{Claim}
\begin{proof}
    For each $b\in C^{1}\sqcup D^{1}\sqcup F^1$,
we can find $(\widehat{x}_{b},\widehat{t}_{b})\in V_{\delta,W,r_{b}}(x_{b},t_{b})\subseteq P^{\ast}(x_{b},t_{b};\Lambda\sqrt{\beta}\lambda,\beta\lambda^2,\beta\lambda^2)$. If $C^{1}\sqcup D^{1}\sqcup F^1=\emptyset$, we are done, so suppose that this is not the case, and fix $(y_{0},s_{0})\in\{(x_{c},t_{c})\}_{c\in C^{1}}\cup\{(x_{d},t_{d})\}_{d\in D^{1}}\cup\{(x_{f},t_{f})\}_{f\in F^{1}}.$ Take some $b\in C^{1}\sqcup D^{1}\sqcup F^1$. Since $V_{\delta,W,1}(x_0,0)$
is not $(\alpha,1)$-temporally independent, we have 
\[
|\widehat{t}_{b}-\widehat{s}_{0}|\leq\alpha.
\]
Therefore, by \cite[Proposition 10.8]{Bam3} and Lemma \ref{lem:weaksplittingepsreg}, we must have that 
$$d_{W_{1}}^{g_{s_{0}-2\lambda^{2}}}(\nu_{\widehat{x}_{b},\widehat{t}_{b};s_{0}-2\lambda^{2}},\nu_{\widehat{x}_{b'},\widehat{t}_{b'};s_{0}-2\lambda^{2}})\leq D_0\lambda.$$
In particular, $(x_b,t_b)\in P^\ast(y_0,s_0;4\Lambda\sqrt{\beta}\lambda,2\lambda^2,2\lambda^2)$. From which it follows from an elementary covering argument that there are most $C(Y,\Lambda)$ elements of $ C^{1}\sqcup D^{1}\sqcup F^1$.
\end{proof}

Using Claim \ref{claim:d_ball}, we get that
\[
\sum_{c\in C^{1}}r_{c}^{2}+\sum_{d\in D^{1}}r_{d}^{2}\leq C(Y,\Lambda)\alpha^{\frac{1}{4}}.
\]
For every  $d\in D^{1}$, $P^{\ast}(x_{d},t_{d};r_{d})$ satisfies the hypotheses of
the lemma, so we may repeat
the construction on each such ball. That is, we obtain a cover
\[
P^{\ast}(x_{d},t_{d};r_{d})\subseteq\bigcup_{c\in C_{d}^{2}}P^{\ast}(x_{c},t_{c};r_{c})\cup\bigcup_{d'\in D_{d}^{2}}P^{\ast}(x_{d},t_{d};r_{d})\cup\bigcup_{e\in E_{d}^{2}}P^{\ast}(x_{e},t_{e};r_{e}),
\]
where $r_{c}=r_{d'}=r_{e}=\left(\frac{\sqrt{\beta}}{2}\lambda\right)^{2}$,
$\sum_{e\in E_{d}^{2}}r_{e}^{2}\leq C_{b}\alpha^{-2}r_{d}^2$, and 
\[
\sum_{c\in C_{d}^{2}}r_{c}^{2}+\sum_{d'\in D_{d}^{2}}r_{d}^{2}\leq C(Y,\Lambda)\alpha^{\frac{1}{4}} r_{d}^{2}.
\]
Setting $D^{2}:=\bigcup_{d\in D^{1}}D_{d}^{2}$,
$C^{2}:=C^{1}\cup\bigcup_{d\in D^{1}}C_{d}^{2}$, $E^{2}:=E^{1}\cup\bigcup_{d\in D^{1}}E_{d}^{2}$,
we obtain a cover
\[
P^{\ast}(x_{\ast},0;1)\subseteq\bigcup_{c\in C^{2}}P^{\ast}(x_{c},t_{d};r_{c})\cup\bigcup_{d\in D^{2}}P^{\ast}(x_{d},t_{d};r_{d})\cup\bigcup_{e\in E^{2}}P^{\ast}(x_{e},t_{e};r_{e}),
\]
where we have $r_{d}=\left(\frac{\sqrt{\beta}}{2}\lambda\right)^{2}$
for all $d\in D^{2}$, as well as the content estimates
\[
\sum_{c\in C^{2}}r_{c}^{2}\leq\sum_{c\in C^{1}}r_{c}^{2}+\sum_{d\in D^{1}}\sum_{c\in C_{d}^{2}}r_{c}^{2}\leq C(Y,\Lambda)\alpha^{\frac{1}{4}}+(C(Y,\Lambda)\alpha^{\frac{1}{4}})^{2},
\]
\[
\sum_{d\in D^{2}}r_{d}^{2}\leq(C(Y,\Lambda)\alpha^{\frac{1}{4}})^{2},
\]
\[
\sum_{e\in E^{2}}r_{e}^{2}\leq C_{b}\alpha^{-2}+(C_{b}\alpha^{-2})(C_{1}\alpha^{\frac{1}{4}}).
\]
Iterating this procedure we immediately get the following.
\begin{Claim} \label{claim:d_ball_3}
    There exists a sequence of coverings
    $$ P^\ast(x_0,0;1)\subseteq \bigcup_{c\in C^k}P^\ast(x_c,t_c;r_c)\cup \bigcup_{d\in D^k}P^\ast(x_d,t_d;r_d)\cup \bigcup_{e\in E^k}P^\ast(x_e,t_e;r_e)\cup\bigcup_{f\in F^k}P^\ast(x_f,t_f;r_f) $$
    where 
    \begin{enumerate}
        \item for every $c\in C^k$, $P^\ast(x_c,t_c;r_c)$ is a $(c)$-ball,
         \item for every $d\in D^k$, $P^\ast(x_d,t_d;r_d)$ is a $(d)$-ball,
          \item for every $e\in E^k$, $P^\ast(x_e,t_e;r_e)$ is an $(e)$-ball,
           \item for every $f\in F^k$, $P^\ast(x_f,t_f;r_f)$ is an $(f)$-ball,
    \end{enumerate}
     $r_d= (\frac{\sqrt{\beta}}{2}\lambda)^k$ for each $d\in D^k$, 
     and we have the content estimates
    \begin{align*}
       \sum_{c\in C^k_d}r_c^2+\sum_{f\in F^k_d}r_f^2&\leq C(Y,\Lambda)\alpha^{\frac{1}{4}}\cdot r_d^2\\ \sum_{e\in E^k_d}r_e^2&\leq C(Y,\alpha)r_d^2\\ \sum_{c\in C^k}r_c^2+\sum_{f\in F^k}r_f^2 &\leq \sum_{l=1}^k(C(Y,\Lambda)\alpha)^{\frac{l}{4}}\\\sum_{e\in E^k}r_e^2& \leq C(Y,\alpha)\sum_{l=0}^{k-1}(C(Y,\Lambda)\alpha)^{\frac{l}{4}}.
    \end{align*}
\end{Claim}
  Since $(M^4,(g_t)_{t\in [-2\delta^{-1},10]})$ is smooth and compact, there is some $k_0\in \mathbb{N}$ for which the $k_0$-th element of the sequence constructed in Claim \ref{claim:d_ball_3} contains no $(d)$-balls. Clearly this covering satisfies all the conditions of the lemma if $\alpha\leq \overline{\alpha}(Y,\Lambda)$, which completes the proof.
\end{proof}


The following Lemma allows us to reduce both the $(f)$-ball, and refined $(c)$-ball, covering to Lemma \ref{lem:real_analysis}.

\begin{Lemma}\label{lem:col_est}
    Let $Y<\infty$ and $\sigma,\gamma,\eta,\alpha>0$. The following statement is true for $\beta\leq \overline{\beta}(Y)$, $\Lambda\geq \underline{\Lambda}(Y)$ and $\delta \leq \overline{\delta}(Y,\sigma,\gamma,\eta,\alpha,\Lambda,\beta)$. Let $(M^4,(g_t)_{t\in [-2\delta^{-1},10]} )$ be a closed K\"ahler Ricci flow and $x_0\in M $ be such that  $\inf_{(x,t)\in P^\ast(x_0,0;2\Lambda)}\mathcal{N}_{(x,t)}(\delta^{-1})\geq -W-\delta $ for some $\sigma \leq W\leq Y-\delta$. If $\{P^\ast(x_i,t_i;r_i)\}_{i\in  F}\subseteq P^\ast(x_0,0;2\Lambda,2,2)$ is any collection of $(f)$-balls such that
    \begin{enumerate}
        \item  for each $f\in F$ we have $r_f=\beta^k$ for some $k\in \mathbb{N}$,
        \item for each $f,f^\prime \in F$ we have $P^\ast(x_f,t_f;\gamma^3r_f)\cap P^\ast(x_{f^\prime},t_{f^\prime};\gamma^3r_{f^\prime})=\emptyset $ if $r_f=r_{f^\prime}$,
        \item for each $f\in F$ we have $|\mathcal{N}_{(x_i,t_i)}(\delta \beta^{-2}r_i^2)+W|\leq \delta$.
    \end{enumerate}
 Then the temportal subsets $\mathfrak{t}(P^\ast(x_f,t_f;r_f))$ together with $$T_r=\overline{\mathfrak{t}\left(\bigcup_{(x,t)\in \bigcup_{i\in F}P^\ast(x_i,t_i;r_i)}V_{\delta,W,r}(x,t)  \right)}\subseteq \mathbb{R}$$ satisfy the  hypothesis of Lemma \ref{lem:real_analysis}.

\end{Lemma}

\begin{proof}
 As in the proof of Lemma \ref{lem:annular}, given $(x_i,t_i),(x_j,t_j)\in  F$, we let $\kappa((x_i,t_i),(x_j,t_j))>0$ be minimal such that $(x_i ,t _i)\in P^{\ast}(x_j,t_j;\kappa)$. Note that  $\kappa((x_i,t_i),(x_j,t_j)) \geq c(\gamma)r_i$ if $r_i=r_j$. By \cite[Proposition 7.1]{Bam3}, for $\delta^\prime>0$ and  $\delta \leq \overline{\delta}(Y,\delta^\prime,\gamma,\beta)$, we have that for $(x_i,t_i)$ is $(\delta^\prime,r)$-selfsimilar for $(x_i,t_i)\in  F$ and $r\in [\kappa_i,1]$.
 We have the following two claims.
\begin{Claim}\label{claim:ball_type_1}
  For $\delta \leq \overline{\delta}(Y,\sigma,\gamma,\eta,\alpha,\Lambda,\beta)$ the following statement is true. For each $(x_i,t_i),(x_j,t_j)\in  F$ with $r_i=r_j$ we have  $|t_i-t_j|\geq c(\gamma,Y,\Lambda )r_i^2$.
\end{Claim}
\begin{proof}
   Suppose that the claim was false, so that after diagonalisation and relabeling there is a sequence $\delta_i\searrow 0$ together with closed K\"ahler Ricci flows $(M^4_i,(g_{i,t})_{t\in [-2\delta^{-1}_i,10]} )$ with $x_{0,i}\in M $ such that  $\inf_{(x,t)\in P^\ast(x_{0,i},0;2\Lambda)}\mathcal{N}_{(x,t)}(\delta_i^{-1})\geq -W-\delta_i $ for some $\sigma \leq W\leq Y-\delta_i$, for which there are $P^\ast$-balls $P^\ast(y_{j,i},s_{j,i};r_{j,i})\subseteq  P^\ast(x_{0,i},0;2\Lambda,2,2)$, $j=1,2$, where $$P^\ast(y_{1,i},s_{1,i};\gamma^3 r_{1,i})\cap P^\ast(y_{2,i},s_{2,i};\gamma^3 r_{2,i})=\emptyset,$$  $(y_{j,i},s_{j,i})$ are $(\delta_i,r)$-selfsimilar for all $r\in [\kappa_i,1]$, $\mathcal{N}_{(y_{j,i},s_{j,i})}(\delta_i\kappa_i^2)\leq -2\sigma$, but  $|t_{1,i}-t_{2,i}|\leq c\cdot r_i^2$. We will derive a contradiction by choosing $c$ to be sufficiently small.  By \cite[Proposition 10.8, Proposition 12.1]{Bam3}, for any $\epsilon>0$, we can choose $c\leq \overline{c}(Y,\gamma,\epsilon,\Lambda)$, so that after rescaling by $\kappa_i$ and relabeling, we have that $(y_{1,i},s_{1,i})$ are strongly $(1,\epsilon,1)$-split. By Lemma \ref{lem:weaksplittingepsreg}, we can choose $\epsilon\leq \overline{\epsilon}(Y)$ so that $\lim_{i\to\infty}\mathcal{N}_{(y_{1,i},s_{1,i})}(\delta_i)=0$, which is a contradiction.
\end{proof}

\begin{Claim}\label{claim:ball_type_2}
    For $\beta\leq \overline{\beta}(Y)$,  $\Lambda\geq \underline{\Lambda}(Y)$, and $\delta\leq \overline{\delta}(Y,\sigma,\gamma,\eta,\alpha,\Lambda,\beta)$, the following statement is true. If $(y_0,s_0)\in \bigcup_{(x,t)\in \bigcup_{i\in F}P^\ast(x_i,t_i;r_i)}V_{\delta,W,r_f}(x,t)$ for some $f\in F$, and $s_0\in [t_f-2r_f^2,t_f+2r_f^2]$, then $(y_0,s_0)\in P^\ast(x_f,t_f;\Lambda r_f,2 r_f^2,2r_f^2)$. In particular, for each $f\in F$ we have $I_f\cap T_{r_f}=\emptyset$. 
\end{Claim}
\begin{proof}
    Take $(y_0,s_0)$ as in the statement of the claim and take $(x_f,t_f)$ such that $(y_0,s_0)\in \bigcup_{(x,t)\in \bigcup_{i\in F}P^\ast(x_i,t_i;r_i)}V_{\delta,W,r_f}(x,t)$. Recall that for $\delta^\prime>0$ both $(y_0,s_0)$ and $(x_f,t_f)$ are $(\delta^\prime,r)$-selfsimilar for each $r\in [ r_f,1]$ if $\delta \leq\overline{\delta}(Y,\beta,\delta^\prime)$. Write $r_f=\beta^k$ for some $k\in\mathbb{N}$. If $k=0$, there is nothing to do, so suppose that $k\geq 1$. Then $(y_0,s_0)\in P^\ast(x_f,t_f;4\Lambda ,2 ,2)$ and $|s_0-t_f|<\beta^2$, and therefore, in a similar fashion to Claim \ref{claim:ball_type_2}, we get from \cite[Proposition 10.8]{Bam3} and Lemma \ref{lem:weaksplittingepsreg} that $(y_0,s_0)\in P^\ast(x_f,t_f;C(Y)\beta,2\beta^2,2\beta^2)$ if $\beta\leq \overline{\beta}(Y)$ and $\delta \leq \overline{\delta}(Y,\sigma,\gamma,\eta,\alpha,\Lambda,\beta)$. Repeating this argument finitely many times, we get that $(y_0,s_0)\in P^\ast(x_f,t_f;C(Y) \beta^k,2\beta^{2k},2\beta^{2k})$. This gives the result if $\Lambda \geq C(Y)$.
\end{proof}

It is easy to see that \eqref{claim:ball_type_1}-\eqref{claim:ball_type_2} imply that the hypothesis of Lemma \ref{lem:real_analysis} are satisfied, and so the proof is complete.
\end{proof}

\begin{Lemma}[$(f)$-ball covering] \label{lem:f_ball}
    Let $\sigma, \alpha,\eta>0$ and $Y<\infty$. For $\Lambda\geq \underline{\Lambda}(Y)$ and  $\delta\leq\overline{\delta }(Y,\sigma,\alpha,\eta,\Lambda)$, the following statement holds.  Let $(M^4,(g_t)_{t\in [-2\delta^{-1},10]} )$ be a closed K\"ahler Ricci flow and $x_0\in M$ be such that $\inf_{P^\ast(x_0,0;2\Lambda)}\mathcal{N}_{(x_0,0)}(\delta^{-1})\geq -W-\delta$ for some $\sigma \leq W\leq Y-\delta$. If $P^\ast(x_0,0;1)$ is an $(f)$-ball, then there exists a covering
    $$ P^\ast(x_0,0;1)\subseteq \bigcup_{c\in C}P^\ast(x_c,t_c;r_c)\cup \bigcup_{e\in E}P^\ast(x_e,t_e;r_e)\cup\bigcup_{d\in D}P^\ast(x_d,t_d;r_d)$$
   where 
    \begin{enumerate}
    \item for every $c\in C$, $P^\ast(x_c,t_c;r_c)$ is a $(c)$-ball, 
       \item for every $d\in D$, $P^\ast(x_d,t_d;r_d)$ is a $(d)$-ball,
         \item for every $e\in E$, $P^\ast(x_e,t_e;r_e)$ is an $(e)$-ball,
    \end{enumerate}
    and we have the content estimates
    \begin{align*}
      \sum_{c\in C}r_c^2+  \sum_{d\in D}r_d^2+\sum_{e\in E}r_e^2& \leq C(Y,\eta,\Lambda). 
    \end{align*}
    
\end{Lemma}

\begin{proof}
Let $\gamma>0$ be a constant to be determined and let $\beta\leq \overline{\beta}(Y)$ be as in Lemma \ref{lem:col_est}. Choose a maximal subset $\{(x_i ,t_i)\}_{i\in \tilde{N}^1}\in V_{\delta,W,1}(x_0,0)\cap P^\ast(x_0,0;\frac{5}{4})$ such that $\{P^\ast(x_i,t_i; \gamma^3\cdot\beta)\}_{i\in\tilde{N}^1}$ are pairwise disjoint, and a maximal subset $\{(x_i ,t_i)\}_{i\in \hat{N}^1}\subseteq P^\ast(x_0,0;1)\backslash T_{10\gamma^3\beta}\left(V_{\delta,W,1}(x_0,0)\right)$ such that 
$\{P^\ast(x_i,t_i; \frac{1}{40\Lambda}\gamma^6\cdot\beta)\}_{i\in\hat{N}^1}$ are pairwise disjoint.  Set $r_i\coloneqq \beta $ for $i\in \tilde{N}^1$ and $r_i=\frac{\gamma^3}{40\Lambda}\beta$ for $i\in \hat{N}^1$. Write $$N^{1}\coloneqq \tilde{N}^1\sqcup \hat{N}^1 =  C^1\sqcup D^1\sqcup E^1\sqcup F^1$$ such that 
\begin{enumerate}
    \item for every $c\in C^1$, $P^\ast(x_c,t_c;r_c)$ is a $(c)$-ball,
      \item for every $d\in D^1$, $P^\ast(x_d,t_d;r_d)$ is a $(d)$-ball,
        \item for every $e\in E^1$, $P^\ast(x_e,t_e;r_e)$ is an $(e)$-ball,
          \item for every $f\in F^1$, $P^\ast(x_f,t_f;r_f)$ is an $(f)$-ball.
\end{enumerate}
Note that by construction, for every $f\in F^1$, $(x_f,t_f)\in V_{\delta,W,1}(x_0,0)$. By an elementary covering argument, we have the content estimate \begin{align*}
	\sum_{f\in F^1}r_f^2+\sum_{c\in C^1}r_c^2+ \sum_{d\in D^1}r_d^2+\sum_{e\in E^1}r_e^2\leq C(Y,\gamma,\Lambda)\beta^{-4}.
\end{align*}

Since every $(f)$-ball $P^\ast(x_f,t_f;r_f)$ satisfies the hypothesis of the lemma, we may continue to decompose every such ball. That is, choose a maximal subset  $ \{(x_i,t_i)\}_{i\in \tilde{N}^2} \subseteq \bigcup_{f\in F^1}V_{\delta,W,\beta}(x_f,t_f)\cap P^\ast(x_f,t_f;\frac{5}{4}r_f)$ such that $\{P^\ast(x_i,t_i; \gamma^3\beta^2)\}_{i\in \tilde{N}^2}$ are pairwise disjoint, and a maximal subset $ \{(x_i,t_i)\}_{i\in \hat{N}^2} \subseteq \bigcup_{f\in F^1}P^\ast(x_f,t_f;r_f) \backslash T_{10\gamma^3\beta^2}\left(\bigcup_{f\in F^1}V_{\delta,W,\beta}(x_f,t_f)\right)$ such that  $\{P^\ast(x_i,t_i; \frac{1}{40\Lambda}\gamma^6\beta^2)\}_{i\in \tilde{N}^2}$ are pairwise disjoint. Set $r_i\coloneqq \beta^2 $ for $i\in \tilde{N}^2$ and $r_i=\frac{\gamma^3}{40\Lambda}\beta^2$ for $i\in \hat{N}^2$. Write $N^2=\tilde{N}^2\sqcup \hat{N}^2= C^2\sqcup D^2\sqcup E^2\sqcup F^2$ so that 
\begin{enumerate}
    \item for every $c\in C^2_f$, $P^\ast(x_c,t_c;r_c)$ is a $(c)$-ball,
      \item for every $d\in D^2_f$, $P^\ast(x_d,t_d;r_d)$ is a $(d)$-ball,
        \item for every $e\in E^2_f$, $P^\ast(x_e,t_e;r_e)$ is an $(e)$-ball,
          \item for every $f\in F^2_f$, $P^\ast(x_f,t_f;r_f)$ is an $(f)$-ball.
\end{enumerate}
In total, this gives a cover 
\begin{align*}
P^\ast(x_0,0;1)\subseteq &\bigcup_{f\in F^1}\left(\bigcup_{c\in C_f^2} P^\ast(x_c,t_c;r_c)\cup \bigcup_{d\in D_f^2} P^\ast(x_d,t_d;r_d)\cup \bigcup_{e\in E_f^2} P^\ast(x_e,t_e;r_e)\cup \bigcup_{f^\prime \in F_f^2} P^\ast(x_{f^\prime},t_{f^\prime};r_f)  \right)\\ & \cup \bigcup_{c\in C^1}P^\ast(x_c,t_c;r_c) \cup \bigcup_{d\in D^1}P^\ast(x_d,t_d;r_d) \cup \bigcup_{e\in E^1}P^\ast(x_e,t_e;r_e)\\&\coloneqq \bigcup_{c\in C^2}P^\ast(x_c,t_c;r_c)\cup \bigcup_{d\in D^2}P^\ast(x_d,t_d;r_d)\cup\bigcup_{e\in E^2}P^\ast(x_e,t_e;r_e)\cup\bigcup_{f\in F^2}P^\ast(x_f,t_f;r_f) ,
\end{align*}
where, by construction, $P^\ast(x_{f},t_f;\gamma^3 r_f)\cap P^\ast(x_{f^\prime},t_{f^\prime};\gamma^3r_{f^\prime})=\emptyset$ if $r_f=r_{f^\prime}$, and $(x_f,t_f)\in V_{\delta,W,\beta^{-1}r_f}(x_{f^\prime},t_{f^\prime})$ for some $(x_{f^\prime},t_{f^\prime})\in F^1$.

Iterating this procedure, we immediately have the following claim.
\begin{Claim} \label{claim:f_ball_cov_1}
	There exists a sequence of coverings
	$$P^\ast(x_0,0;1)\subseteq \bigcup_{c\in C^k}P^\ast(x_c,t_c;r_c)\cup  \bigcup_{d\in D^k}P^\ast(x_d,t_d;r_d)\cup  \bigcup_{e\in E^k}P^\ast(x_e,t_e;r_e)\cup  \bigcup_{f\in F^k}P^\ast(x_f,t_f;r_f) $$
	 where $r_f=\beta^k$ for each $f\in F^k$,
	\begin{align}
		\sum_{f\in F^1}r_f^2+\sum_{c\in C^1}r_c^2+ \sum_{d\in D^1}r_d^2+\sum_{e\in E^1}r_e^2&\leq C(Y,\gamma,\Lambda)\beta^{-4}, \label{iden:f_ball_cov_3}
	\end{align}
   and 
   \begin{enumerate}
       \item for each $f,f^\prime \in \bigcup_{l=1}^k\bigcup_{\hat{f}\in F^l}F^l$ we have  $P^\ast(x_{f},t_f;\gamma^3 r_f)\cap P^\ast(x_{f^\prime},t_{f^\prime};\gamma^3r_{f^\prime})=\emptyset$ if $r_f=r_{f^\prime}$,
       \item  for every $f\in F^k$ we have $(x_f,t_f)\in V_{\delta,W,\beta^{-1}r_f}(x_{f^\prime},t_{f^\prime})$ for some $(x_{f^\prime},t_{f^\prime})\in F^{k-1}$.
   \end{enumerate}
\end{Claim}
 Consider the sequence of coverings as in Claim \ref{claim:f_ball_cov_1}.  
Combining \eqref{iden:f_ball_cov_3} with Lemma \ref{lem:col_est} and Lemma \ref{lem:real_analysis} gives
  \begin{align*}
       \sum_{l=1}^{k}\sum_{c\in C^l}r_c^2+ \sum_{l=1}^{k}\sum_{d\in D^l}r_d^2+  \sum_{l=1}^{k}\sum_{e\in E^l}r_e^2&\leq \left(C(\gamma,Y,\beta)\sum_{l=1}^{k-1}\sum_{f\in F^l}r_f^2+\sum_{c\in C^1}r_c^2+\sum_{d\in D^1}r_d^2+\sum_{e\in E^1}r_e^2\right)\\& \leq C(\eta,\gamma,Y,\beta,\Lambda)
  \end{align*} 
 if $\Lambda\geq \underline{\Lambda}(Y)$ and $\delta\leq\overline{\delta }(Y,\sigma,\alpha,\eta,\Lambda,\beta)$.   Since $(M^4,(g_t)_{t\in [-2\delta^{-1},10]})$ is smooth and compact, there is some $k_0\in \mathbb{N}$ for which the $k_0$-th element of the sequence constructed in Claim \ref{claim:f_ball_cov_1} contains no $(f)$-balls. We conclude that this covering satisfies all the conditions of the lemma, which completes the proof.
\end{proof}

Next, we combine Lemma \ref{lem:f_ball} with Lemma \ref{lem:coarse_c} and  Lemma \ref{lem:coarse_d} to produce refined $(c)$ and $(d)$-ball covers that do not contain any $(f)$-balls. 

\begin{Lemma}[Refined $(c)$-ball covering]\label{lem:ref_c}
      Let $\sigma ,\alpha>0$ and  $Y<\infty$. Then for $\Lambda\geq \underline{\Lambda}(Y)$, $\gamma\leq \overline{\gamma}(Y,\Lambda)$, $\epsilon\leq \overline{\epsilon}(Y,\gamma,\Lambda)$, $\eta\leq \overline{\eta}(Y,\gamma,\Lambda)$, and $\delta \leq \overline{\delta}(Y,\sigma,\epsilon,\gamma,\eta,\alpha,\Lambda)$, the following statement is true.  Let $(M^4,(g_t)_{t\in [-2\delta^{-1},10]} )$ be a closed K\"ahler Ricci flow and $x_0\in M $ be such that  $\inf_{(x,t)\in P^\ast(x_0,0;2\Lambda)}\mathcal{N}_{(x,t)}(\delta^{-1})\geq -W-\delta $ for some $\sigma \leq W\leq Y-\delta$. Then if $P^\ast(x_0,0;1)$ is a $(c)$-ball, there exists a covering
    $$ P^\ast(x_0,t_0;1)\subseteq \bigcup_{i\in N}\mathcal{A}^i\cup \bigcup_{d\in D}P^\ast(x_d,t_d;r_d)\cup \bigcup_{e\in E}P^\ast(x_e,t_e;r_e)$$
    where 
   \begin{enumerate}
       \item for every $i\in N$, $\mathcal{A}^i=P^\ast(x_i,t_i;2\Lambda r_i,2r_i^2,2r_i^2)\backslash\bigcup_{j\in\mathcal{C}^i}P^\ast(x_j,t_j;r_j)$ is an $(m,\epsilon,\gamma,\Lambda)$-annular region,
       \item for every $d\in D$, $P^\ast(x_d,t_d;r_d)$ is a $(d)$-ball,
       \item for every $e\in E$, $P^\ast(x_e,t_e;r_e)$ is an $(e)$-ball, 
   \end{enumerate}
   and we have the content estimate
   \begin{align*}
       \sum_{i\in N}r_i^2+\sum_{d\in D}r_d^2+\sum_{e\in E}r_e^2\leq C(Y,\gamma,\eta,\Lambda). 
   \end{align*}
\end{Lemma}
\begin{proof}
   Let $\beta\leq \overline{\beta}(Y)$ be as in Lemma \ref{lem:col_est}. Since $P^\ast(x_0,0;1)$ is a $(c)$-ball,  for  $\eta\leq \overline{\eta}(Y,\gamma,\beta)$ and $\delta \leq \overline{\delta}(Y,\epsilon,\gamma,\alpha,\eta,\Lambda,\beta)$, we may apply Lemma \ref{lem:coarse_c} (with $\tilde{F}=\emptyset$) to get a covering
    $$P^\ast(x_0,0;1)\subseteq \mathcal{A}^1\cup \bigcup_{d\in D^1}P^\ast(x_d,t_d;r_d)\cup \bigcup_{e\in E^1}P^\ast(x_e,t_e;r_e)\cup \bigcup_{f\in F^1}P^\ast(x_f,t_f;r_f)$$
    where 
    \begin{enumerate}
    \item $\mathcal{A}^1=P^\ast(x_0,0;2\Lambda,2,2)\backslash \bigcup_{i\in \mathcal{C}^1}P^\ast(x_i,t_i;r_i)$, $\mathcal{C}^1= D^1\cup E^1\cup F^1$, is an $(m,\epsilon,\gamma,\Lambda)$-annular region,
      \item for every $d\in D^1$, $P^\ast(x_d,t_d;r_d)$ is a $(d)$-ball,
        \item for every $e\in E^1$, $P^\ast(x_e,t_e;r_e)$ is an $(e)$-ball, 
          \item for every $f\in F^1$, $P^\ast(x_f,t_f;r_f)$ is an $(f)$-ball.
    \end{enumerate}
     By \eqref{iden:annular_struc_1} of Lemma \ref{lem:annular}, we have 
     \begin{align}
     	\sum_{d\in D^1}r_d^2+\sum_{e\in E^1}r_e^2+\sum_{f\in F^1}r_f^2\leq C(Y,\gamma) \label{iden:ref_c_proof_1}
     \end{align}
  if $\gamma \leq \overline{\gamma}(Y,\Lambda)$ and $\epsilon\leq \overline{\epsilon}(Y,\gamma,\Lambda)$.   By Lemma \ref{lem:f_ball}, we can choose $\Lambda \geq \underline{\Lambda}(Y)$ and $\delta \leq \overline{\delta}(Y,\eta,\alpha,\sigma,\Lambda)$ so that for every $f\in F^1$, there exists a cover
    $$P^\ast(x_f,f_f;r_f)\subseteq \bigcup_{c\in C_f^2}P^\ast(x_c,t_c;r_c)\cup \bigcup_{e\in E^2_f}P^\ast(x_e,t_e;r_e)\cup\bigcup_{d\in D^2_f}P^\ast(x_d,t_d;r_d)$$
    where 
       \begin{enumerate}
      \item for every $c\in C^2_f$, $P^\ast(x_c,t_c;r_c)$ is a $(c)$-ball,
          \item for every $d\in D^2_f$, $P^\ast(x_d,t_d;r_d)$ is a $(d)$-ball,
                \item for every $e\in E^2_f$, $P^\ast(x_e,t_e;r_e)$ is a $(e)$-ball,
    \end{enumerate}
   and we have the content estimates  \begin{align}
      \sum_{d\in D^2_f}r_d^2+\sum_{e\in E^2_f}r_e^2+  \sum_{c\in C^2_f}r_c^2&\leq C(Y,\eta,\Lambda)r_f^2. \label{iden:ref_c_proof_2}
    \end{align} 
 Write $C^2=\bigcup_{f\in F^1}C^2_f=\{c_i\}_{i=0}^N$ for some $N\in \mathbb{N}$. Then $P^\ast(x_{c_0},t_{c_0};r_{c_0})$ is a $(c)$-ball and so we may apply Lemma \ref{lem:coarse_c} with $\tilde{F}_{c_0}^2=F^1$ to get a covering
    $$P^\ast(x_{c_0},t_{c_0};r_{c_0})\subseteq \mathcal{A}^2_{c_0}\cup \bigcup_{d\in D^2_{c_0}}P^\ast(x_d,t_d;r_d)\cup \bigcup_{e\in E^2_{c_0}}P^\ast(x_e,t_e;r_e)\cup \bigcup_{f\in F^2_{c_0}}P^\ast(x_f,t_f;r_f)$$
    where
       \begin{enumerate}
       \item $\mathcal{A}^2_{c_0}=P^\ast(x_{c_0},t_{c_0};2\Lambda r_{c_0},2r_{c_0}^2,2r_{c_0}^2)\backslash \bigcup _{i\in \mathcal{C}^2_c}P^\ast(x_i,t_i;r_i)$, $\mathcal{C}^2_{c_0}=D^2_{c_0}\cup E^2_{c_0}\cup F^2_{c_0}$, is an $(m,\epsilon,\gamma,\Lambda)$-annular region,
          \item for every $d\in D^2_{c_0}$, $P^\ast(x_d,t_d;r_d)$ is a $(d)$-ball,
                \item for every $e\in E^2_{c_0}$, $P^\ast(x_e,t_e;r_e)$ is an $(e)$-ball,
                  \item for every $f\in F^2_{c_0}$, $P^\ast(x_f,t_f;r_f)$ is an $(f)$-ball.
    \end{enumerate}
    Next, since $P^\ast(x_{c_1},t_{c_1};r_{c_1})$ is a $(c)$-ball,  we may apply Lemma \ref{lem:coarse_c} with $\tilde{F}^2_{c_1}=F^1\cup F^2_{c_0}$ to get a covering
    $$P^\ast(x_{c_1},t_{c_1};r_{c_1})\subseteq \mathcal{A}^2_{c_1}\cup \bigcup_{d\in D^2_{c_1}}P^\ast(x_d,t_d;r_d)\cup \bigcup_{e\in E^2_{c_1}}P^\ast(x_e,t_e;r_e)\cup \bigcup_{f\in F^2_{c_1}}P^\ast(x_f,t_f;r_f)$$
    where
       \begin{enumerate}
       \item $\mathcal{A}^2_{c_1}=P^\ast(x_{c_1},t_{c_1};2\Lambda r_{c_1},2r_{c_1}^2,2r_{c_1}^2)\backslash \bigcup _{i\in \mathcal{C}^2_{c_1}}P^\ast(x_i,t_i;r_i)$, $\mathcal{C}^2_{c_1}=D^2_{c_1}\cup E^2_{c_1}\cup F^2_{c_1}$, is an $(m,\epsilon,\gamma)$-annular region,
          \item for every $d\in D^2_{c_1}$, $P^\ast(x_d,t_d;r_d)$ is a $(d)$-ball,
                \item for every $e\in E^2_{c_1}$, $P^\ast(x_e,t_e;r_e)$ is an $(e)$-ball,
                  \item for every $f\in F^2_{c_1}$, $P^\ast(x_f,t_f;r_f)$ is an $(f)$-ball.
    \end{enumerate}
    Repeating this construction over $C^2$, for $\Lambda\geq \underline{\Lambda}(Y)$ and  $\delta\leq\overline{\delta}(Y,\epsilon,\alpha,\gamma,\eta,\Lambda)$ we now have the following covering of our original ball
    \begin{align*} P^\ast(x_0,0;1)\subseteq & \mathcal{A}^1\cup \bigcup_{d\in D^1}P^\ast(x_d,t_d;r_d)\cup \bigcup_{e\in E^1}P^\ast(x_e,t_e;r_e)\\  \cup& \bigcup_{f\in F^1}\left( \bigcup_{e\in E^2_f}P^\ast(x_e,t_e;r_e)\cup\bigcup_{d\in D^2_f}P^\ast(x_d,t_d;r_d)\right)\\\cup&\bigcup_{f\in F^1}\bigcup_{c\in C^2_f}\left(\mathcal{A}^2_c\cup \bigcup_{d\in D^2_c}P^\ast(x_d,t_d;r_d)\cup \bigcup_{e\in E^2_c}P^\ast(x_e,t_e;r_e)\cup \bigcup_{f\in F^2_c\backslash \tilde{F}^2_c}P^\ast(x_f,t_f;r_f)\right) \\ \coloneqq  &\bigcup_{i\in N^2}\mathcal{A}^i \cup \bigcup _{d\in D^2}P^\ast(x_d,t_d;r_d)\cup \bigcup_{e\in E^2}P^\ast(x_e,t_e;r_e)\cup \bigcup_{f\in F^2}P^\ast(x_f,t_f;r_f).\end{align*}
   Note that by construction, 
   \begin{enumerate}
       \item  for any $f\in F^1\cup F^2$ we have $r_f=\beta^k$ for some $k\in\mathbb{N}$,
       \item  for any $f\in F^1\cup F^2$ we have $(x_f,t_f)\in V_{\delta,W,\beta^{-1}r_f}(x,t)$ for some $(x,t)\in P^\ast(x_0,0;2\Lambda)$,
       \item for any $f,f^\prime\in F^1\cup F^2$ we have $P^\ast(x_f,t_f;r_f)\cap P^\ast(x_{f^\prime},t_{f^\prime};r_{f^\prime})=\emptyset$ if $r_f=r_{f^\prime}$.
   \end{enumerate}
Iterating this procedure, we have the following claim.
\begin{Claim}\label{claim:ref_c_ball_1}
There is a sequence of covers
$$ P^\ast(x_0,t_0;1)=\bigcup_{i\in N^l}\mathcal{A}^i \cup \bigcup _{d\in D^l}P^\ast(x_d,t_d;r_d)\cup \bigcup_{e\in E^l}P^\ast(x_e,t_e;r_e)\cup \bigcup_{f\in F^l}P^\ast(x_f,t_f;r_f)$$
where 
\begin{enumerate}
	\item for every $i\in N^l$, $\mathcal{A}^i=P^\ast(x_i,t_i;2r_i)\backslash \bigcup_{j\in \mathcal{C}^i}P^\ast(x_j,t_j;r_j)$ is an $(m,\epsilon,\gamma,\Lambda)$-annular region,
	\item for every $d\in D^l$, $P^\ast(x_d,t_d;r_d)$ is a $(d)$-ball,
	\item for every $e\in E^l$, $P^\ast(x_e,t_e;r_e)$ is an $(e)$-ball,
	\item for every $f\in F^l$, $P^\ast(x_f,t_f;r_f)$ is an $(f)$-ball,
\end{enumerate}
$r_f\leq \beta ^k$ for every $f\in F^k$, and 
\begin{enumerate}
    \item  for any $f\in \bigcup_{l=1}^kF^l$ we have $r_f=\beta^j$ for some $j\in\mathbb{N}$,
    \item for any $f\in \bigcup_{l=1}^kF^l$ we have $(x_f,t_f)\in V_{\delta,W,\beta^{-1}r_f}(x,t)$ for some $(x,t)\in P^\ast(x_0,0;2\Lambda)$,
    \item for any $f,f^\prime\in\bigcup_{l=1}^kF^l$ we have $P^\ast(x_f,t_f;\gamma^3r_f)\cap P^\ast(x_{f^\prime},t_{f^\prime};\gamma^3 r_{f^\prime})=\emptyset$ if $r_f=r_{f^\prime}$.
\end{enumerate}
\end{Claim}
Consider the sequence of coverings constructed in Claim \ref{claim:ref_c_ball_1}. By Lemma \ref{lem:col_est} and Lemma \ref{lem:real_analysis}, we have
\begin{align}
 	\sum_{i\in N^k}r_i^2+ \sum_{e\in E^k}r_e^2+\sum_{d\in D^k} r_d^2\leq C(Y,\gamma,\eta,\Lambda)\left(1+\sum_{i=1}^k\sum_{f\in F^i}r_f^2\right)\leq C(Y,\gamma,\Lambda,\eta). \label{iden:ref_c_ball_3}
 \end{align}


Since the $(M,(g_t)_{t\in [-2\delta^{-1},10]})$ is smooth and compact, there is some $k_0\in \mathbb{N}$ such that the $k_0$-th element of the sequence  constructed in Claim \ref{claim:ref_c_ball_1} contains no $(f)$-balls. By \eqref{iden:ref_c_ball_3}, we see that this covering satisfies the requirements of the  lemma, which completes the proof.

\end{proof}

\begin{Lemma}[Refined $(d)$-ball covering] \label{lem:ref_d}
    Let $Y<\infty$ and $\eta,\sigma>0$. Then for $\Lambda\geq \underline{\Lambda}(Y)$, $\alpha\leq \overline{\alpha}(Y,\eta,\Lambda)$, and $\delta \leq \overline{\delta}(Y,\alpha,\sigma,\eta,\Lambda)$ the following statement holds. Let $(M^4,(g_{t})_{t\in [-2\delta^{-1},10]} )$ be a closed K\"ahler Ricci flow, and $x_0\in M$ be such that $\inf_{P^\ast(x_0,0;2\Lambda)}\mathcal{N}_{x_0,t_0}(\delta^{-1})\geq -W-\delta$ for some $\sigma \leq W\leq Y-\delta$. Then if $P^\ast(x_0,0;1)$ is a $(d)$-ball, there exists a covering
    $$ P^\ast(x_0,0;1)\subseteq \bigcup_{c\in C}P^\ast(x_c,t_c;r_c)\cup \bigcup_{e\in E}P^\ast(x_e,t_e;r_e)$$
    where 
    \begin{enumerate}
        \item for every $c\in C$, $P^\ast(x_c,t_c;r_c)$ is a $(c)$-ball,
        \item for every $e\in E$, $P^\ast(x_e,t_e;r_e)$ is an $(e)$-ball,
    \end{enumerate}
    and we have the content estimates
    \begin{align*}
        \sum_{c\in C}r_c^2& \leq C(Y,\Lambda)\alpha^{\frac{1}{4}} \\ \sum_{e\in E}r_e^2& \leq C(Y,\alpha,\Lambda). 
    \end{align*}
\end{Lemma}

\begin{proof}
 By Lemma \ref{lem:coarse_d}, for $\alpha\leq \overline{\alpha}(Y,\Lambda)$ and $\delta\leq \overline{\delta}(\alpha,Y,\sigma,\eta,\Lambda)$ we can produce a cover 
    $$ P^\ast(x_0,0;1)\subseteq  \bigcup_{c\in C^1}P^\ast(x_c,t_c;r_c)\cup \bigcup_{e\in E^1}P^\ast(x_e,t_e;r_e) \cup \bigcup_{f\in F^1}P^\ast(x_f,t_f;r_f)$$
    where 
    \begin{enumerate}
          \item for every $c\in C^1$, $P^\ast(x_c,t_c;r_c)$ is a $(c)$-ball,
        \item for every $d\in D^1$, $P^\ast(x_d,t_d;r_d)$ is a $(d)$-ball,
        \item for every $f\in F^1$, $P^\ast(x_f,t_f;r_f)$ is an $(f)$-ball,
    \end{enumerate}
    and we have the content estimates 
    \begin{align}
        \sum_{c\in C^1}r_c^2+\sum_{f\in F^1}r_f^2&\leq C(Y,\Lambda)\alpha^{\frac{1}{4}} \label{iden:ref_d_proof_1}  \\ \sum_{e\in E^1}r_e^2&\leq C(Y,\alpha,\Lambda). \label{iden:ref_d_proof_2}
    \end{align}
Choose $\delta \leq \overline{\delta}(Y,\eta,\alpha,\Lambda,\sigma)$ so that for every $f\in F^1$,  we may apply Lemma \ref{lem:f_ball} to produce a cover
    $$ P^\ast(x_f,t_f;r_f)\subseteq  \bigcup_{c\in C^2_f}P^\ast(x_c,t_c;r_c)\cup \bigcup_{e\in E^2_f}P^\ast(x_e,t_e;r_e) \cup \bigcup_{d\in D^2_f}P^\ast(x_d,t_d;r_d)$$
    where we have the content estimate \begin{align}\sum_{c\in C^2_f}r_c^2+\sum_{d\in D^2_f}r_d^2+\sum_{e\in E^2_f}r_e^2 & \leq C(Y,\eta,\Lambda)r_f^2.  \label{iden:ref_d_proof_3}\end{align}
    For every $f\in F^1$ and $d\in D^2_f$, we can apply Lemma \ref{lem:coarse_d} once again to produce a cover
    $$ P^\ast(x_d,t_d;r_d)\subseteq  \bigcup_{c\in C^2_d}P^\ast(x_c,t_c;r_c)\cup \bigcup_{e\in E^2_d}P^\ast(x_e,t_e;r_e) \cup \bigcup_{f\in F^2_d}P^\ast(x_f,t_f;r_f)$$
    where we have the content estimates
    \begin{align}
            \sum_{c\in C^2_d}r_c^2+\sum_{f\in F^2_d}r_f^2&\leq C(Y,\Lambda)\alpha^{\frac{1}{4}}\cdot r_d^2 \label{iden:ref_d_proof_4}  \\ \sum_{e\in E^2_d}r_e^2&\leq C(Y,\alpha,\Lambda)r_d^2. \label{iden:ref_d_proof_5}
    \end{align}
    Overall, we have produced a cover
    \begin{align*} P^\ast(x_0,0;1)\subseteq &  \bigcup_{c\in C^1}P^\ast(x_c,t_c;r_c)\cup \bigcup_{e\in E^1}P^\ast(x_e,t_e;r_e) \\ \cup & \bigcup_{f\in F^1}\left( \bigcup_{c\in C^2_f}P^\ast(x_c,t_c;r_c)\cup \bigcup_{e\in E^2_f}P^\ast(x_e,t_e;r_e)\right) \\ \cup & \bigcup_{f\in F^1}\bigcup_{d\in D^2_f}\left( \bigcup_{c\in C^2_d}P^\ast(x_c,t_c;r_c)\cup \bigcup_{e\in E^2_d}P^\ast(x_e,t_e;r_e) \cup \bigcup_{f\in F^2_d}P^\ast(x_f,t_f;r_f)\right) \\  \coloneqq  & \bigcup_{c\in C^2}P^\ast(x_c,t_c;r_c)\cup \bigcup_{e\in E^2}P^\ast(x_e,t_e;r_e)\cup \bigcup_{c\in F^2}P^\ast(x_f,t_f;r_f). \end{align*}
    By \eqref{iden:ref_d_proof_1}-\eqref{iden:ref_d_proof_5} we have
    \begin{align*}
        \sum_{c\in C^2}r_c^2&=\sum_{c\in C^1}r_c^2+\sum_{f\in F^1}\sum_{c\in C^2_f}r_c^2+\sum_{f\in F^1}\sum_{d\in D^2_f}\sum_{c\in C^2_d}r_c^2\\ & \leq C(Y,\eta,\Lambda)\alpha^{\frac{1}{4}} +(C(Y,\eta,\Lambda)\alpha^{\frac{1}{4}})^2\\  \sum_{e\in E^2}&= \sum_{e\in E^1}r_e^2+\sum_{f\in F^1}\sum_{e\in E^2_f}r_e^2+\sum_{f\in F^1}\sum_{d\in D^2_f}\sum_{e\in E^2_d}r_e^2\\ & \leq C(Y,\eta,\Lambda)(1+\alpha^{\frac{1}{4}}).
    \end{align*}
    Iterating the above procedure, we immediately get the following claim.
    \begin{Claim}\label{claim:ref_d}
    There is a sequence of coverings
     $$ P^\ast(x_0,0;1)\subseteq \bigcup_{c\in C^k}P^\ast(x_c,t_c;r_c)\cup \bigcup_{e\in E^k}P^\ast(x_e,t_e;r_e)\cup \bigcup_{c\in F^k}P^\ast(x_f,t_f;r_f) $$
        where 
        \begin{enumerate}
        \item for every $c\in C^k$, $P^\ast(x_c,t_c;r_c)$ is a $(c)$-ball,
        \item for every $e\in E^k$, $P^\ast(x_e,t_e;r_e)$ is an $(e)$-ball,
        \item for every $f\in F^k$, $P^\ast(x_f,t_f;r_f)$ is an $(f)$-ball,
        \end{enumerate}
       $r_f\leq \beta ^{k}$ for every $f\in F^k$,
        such that we have the content estimates
        \begin{align*}
              \sum_{c\in C^k}r_c^2&\leq \sum_{l=1}^k(C(Y,\eta,\Lambda)\alpha^{\frac{1}{4}})^l\\ \sum_{e\in E^k}r_e^2&\leq \sum_{l=0}^{k-1}(C(Y,\eta,\Lambda)\alpha^{\frac{1}{4}})^l.
        \end{align*}
            \end{Claim}
  Since $(M,(g_t)_{t\in [-2\delta^{-1},0]})$ is smooth and compact, there is some $k_0\in \mathbb{N}$ such that the $k_0$-th element of the sequence constructed in Claim \ref{claim:ref_d}  contains no $(f)$-balls. Clearly, this covering satisfies the requirements of the lemma if $\alpha\leq \overline{\alpha}(\eta,Y,\Lambda)$, which completes the proof.


    
 \end{proof}
Next we combine Lemma \ref{lem:ref_c}, Lemma \ref{lem:f_ball}, and Lemma \ref{lem:ref_d} to prove Theorem \ref{thm:neck_decomp}.


\begin{proof}
Let $\sigma>0$ be to be determined. Let $\Lambda\geq \underline{\Lambda}(Y)$ and $\gamma\leq \overline{\gamma}(Y,\Lambda)$ be as in Lemma \ref{lem:ref_c}. Note that we may assume without loss of generality that  $\epsilon\leq \overline{\epsilon}(Y,\Lambda,\gamma)$, also as in Lemma \ref{lem:ref_c}. Choose  $\eta\leq \overline{\eta}(Y,\gamma,\Lambda)$, $\alpha\leq \overline{\alpha}(Y,\eta,\Lambda)$, and  $\delta\coloneqq \overline{\delta}(Y,\epsilon,\alpha,\eta,\sigma,\Lambda,\gamma )$, as in Lemma \ref{lem:col_est}, Lemma \ref{lem:f_ball}, Lemma \ref{lem:ref_c}, and Lemma \ref{lem:ref_d}. By adjusting $Y$ if necessary, we may assume that $\inf_{(x,t)\in P^\ast(x_0,0;2\Lambda)}\geq -Y$. Consider a maximal subset $\{(x_i,t_i)\}_{i\in N^1}\in P^\ast(x_0,0;1)$ such that
 $\{P^\ast(x_i,t_i;\delta^2)\}_{i\in N^1}$ are pairwise disjoint. For $W=Y$, write $N^1=C^1\sqcup D^1\sqcup E^1\sqcup F^1$ so that 
\begin{enumerate}
    \item for every $c\in C^1$, $P^\ast(x_c,t_c;\delta)$ is a $(c)$-ball,
    \item for every $d\in D^1$, $P^\ast(x_d,t_d;\delta)$ is a $(d)$-ball,
     \item for every $e\in E^1$, $P^\ast(x_e,t_e;\delta)$ is an $(e)$-ball,
         \item for every $f\in F^1$, $P^\ast(x_f,t_f;\delta)$ is an $(f)$-ball.
\end{enumerate}
For each $f\in F^1$, we may apply Lemma \ref{lem:f_ball} to get a covering
$$ P^\ast(x_f,t_f;\delta)\subseteq \bigcup_{c\in C^1_f}P^\ast(x_c,t_c;r_c) \cup  \bigcup_{d\in D^1_f}P^\ast(x_d,t_d;r_d)\cup\bigcup_{e\in E^1_f}P^\ast(x_e,t_e;r_e) $$
where we have the content estimates \begin{align}
  \sum_{c\in C^1_f}r_c^2+\sum_{d\in D^1_f}r_d^2+\sum_{e\in E^1_f}r_e^2& \leq C(Y,\eta,\Lambda)\delta^2. \label{iden:neck_decom_proof_1}
\end{align}
Slightly abusing notation we set 
$$\begin{cases}
	C^1=C^1\cup \bigcup_{f\in F^1}C^1_f\\
D^1=D^1\cup \bigcup_{f\in F^1}D^1_f\\
E^1=E^1\cup \bigcup_{f\in F^1}E^1_f
\end{cases},$$
   so that we have a cover
   \begin{align}
P^\ast(x_0,0;1)\subseteq \bigcup_{c\in C^1 }P^\ast(x_c,t_c;r_c)\cup\bigcup_{d\in D^1 }P^\ast(x_d,t_d;r_d)\cup \bigcup_{e\in E^1 }P^\ast(x_e,t_e;r_e) .\label{iden:neck_proof_cov_1}
   \end{align}
Combining \eqref{iden:neck_decom_proof_1} with an elementary covering argument then gives \begin{align}
	\sum _{c\in C^1}r_c^2+\sum_{d\in D^1}r_d^2+\sum_{e\in E^1}r_e^2 \leq C(Y,\eta,\delta).\label{iden:neck_decomp_proof_2}
\end{align} 
We will now remove the $(c)$ and $(d)$-balls from this cover. Take $c\in C^1$, and apply Lemma \ref{lem:ref_c} to get a cover $$ P^\ast(x_c,t_c;r_c)\subseteq \bigcup_{i\in \tilde{N}^1_c}\mathcal{A}_c^i\cup \bigcup_{d\in \tilde{D}_c^1}P^\ast(x_d,t_d;r_d)\cup \bigcup _{e\in \tilde{E}^1_c}P^\ast(x_e,t_e;r_e)$$
where 
\begin{enumerate}
\item for each $i\in \tilde{N}^1_c$, $\mathcal{A}^i_c=P^\ast(x_i,t_i;2\Lambda r_i,2r_i,2r_i)\backslash_{j \in \mathcal{C}^i}P^\ast(x_j,t_j;r_j)$ is an $(m,\epsilon,\gamma,\Lambda)$-annular region for some $m\in \mathbb{N}$ with $m\leq e^Y$,
    \item for each $d\in \tilde{D}^1_c$, $P^\ast(x_d,t_d;r_d)$ is a $(d)$-ball,
    \item  for each $e\in \tilde{E}^1_c$, $P^\ast(x_e,t_e;r_e)$ is an $(e)$-ball,
\end{enumerate}
and  
\begin{align}
	\sum_{i\in \tilde{N}^1_c}r_i^2+\sum_{d\in \tilde{D}^1_c}r_d^2+\sum_{e\in \tilde{E}^1_c}r_e^2\leq C(Y,\eta,\gamma,\Lambda)r_c^2. \label{iden:ref_c_neck_proof}
\end{align}
For every $d\in \tilde{D}^1_c$, we may apply Lemma \ref{lem:ref_d} to get a covering $$P^\ast(x_d,t_d;r_d)\subseteq \bigcup_{c\in \tilde{C}^2_d}P^\ast(x_c,t_c;r_c)\cup \bigcup_{e\in \tilde{E}^2_d}P^\ast(x_e,t_e;r_e),$$
where we have
\begin{align}
	\sum_{c\in \tilde{C}^2_d}r_c^2&\leq C(Y,\eta,\gamma,\Lambda)\alpha^{\frac{1}{4}}r_d^2 \label{iden:ref_c_neck_proof_2}\\
	\sum_{e\in \tilde{E}^2_d}r_e^2&\leq C(Y,\alpha,\eta,\Lambda)r_d^2\label{iden:ref_d_stat_2}.
\end{align}
Setting $\tilde{N}^2_c=\tilde{N}^1_c, \tilde{C}^2_c=\bigcup_{d\in \tilde{D}^1_c}\tilde{C}^2_d,$ and $ \tilde{E}^2_c=E^1_c\cup  \bigcup_{d\in \tilde{D}^1_c}\tilde{E}^2_d$, we now have a cover
$$P^\ast(x_c,t_c;r_c)\subseteq \bigcup_{i\in \tilde{N}^2_c}\mathcal{A}^i\cup \bigcup_{c\in \tilde{C}^2_c}P^\ast(x_c,t_c;r_c)\cup \bigcup_{e\in \tilde{E}^2_c}P^\ast(x_e,t_e;r_e).$$
Combining \eqref{iden:ref_c_neck_proof} and \eqref{iden:ref_d_stat_2} we have
\begin{align*}
 \sum_{i\in \tilde{N}^2_c}r_i^2+  \sum_{e\in \tilde{E}^2_c}r_e^2=\sum_{e\in \tilde{E}^1_c}r_e^2+\sum_{d\in \tilde{D}^1_c}\sum_{e\in \tilde{E}^2_d}r_e^2&\leq C(Y,\gamma,\eta,\Lambda,\alpha)r_c^2\\ \sum_{c^\prime\in \tilde{C}^2_c}r_{c^\prime}^2&\leq C(Y,\gamma,\eta,\Lambda)\alpha^{\frac{1}{4}}\cdot r_c^2.\end{align*}
It is an immediate consequence of the proof of Lemmas \ref{lem:coarse_c}, \ref{lem:coarse_d}, and \ref{lem:f_ball}, that we may assume that $r_{c^\prime}\leq \frac{1}{4}r_c$ for every $c^\prime \in \tilde{C}^2_c$.  In the same way as in the proof of Lemma \ref{lem:coarse_d} and Lemma \ref{lem:ref_d}, if     for $\alpha \leq\overline{\alpha}(Y,\eta,\gamma,\Lambda)$, we iterate this construction to arrive at the following claim. 

\begin{Claim}\label{claim:neck_proof_1}
	For each $c\in C^1$, there is a sequence of coverings 
$$P^\ast(x_c,t_c;r_c)\subseteq \bigcup_{i\in \tilde{N}^k_c}\mathcal{A}^i\cup \bigcup_{c\in \tilde{C}^k_c}P^\ast(x_c,t_c;r_c) \cup \bigcup_{e\in \tilde{E}^k_c}P^\ast(x_e,t_e;r_e)$$
where $r_{c^\prime}\leq \frac{1}{2^{2k}}r_c$ for $c^\prime \in \tilde{C}^k$,
\begin{enumerate}
	\item for every $i\in \tilde{N}^k_c$, $\mathcal{A}^i$ is an $(m,\epsilon,\gamma,\Lambda)$-annular region,
	   \item for every $c^\prime \in \tilde{C}^k_c$, $P^\ast(x_{c^\prime},t_{c^\prime};r_{c^\prime})$ is a $(c)$-ball,
\item for every $e\in \tilde{E}^k_c$, $P^\ast(x_e,t_e;r_e)$ is an $(e)$-ball,
and we have the content estimate 
$$  \sum_{i\in \tilde{N}^k_c}r_i^2+  \sum_{e\in \tilde{E}^k_c}r_e^2\leq \sum_{l=0}^{k-1}(C(Y,\gamma,\eta,\Lambda)\alpha^{\frac{1}{4}})^lr_c^2\leq C(Y,\gamma,\eta,\alpha,\Lambda)r_c^2.$$
\end{enumerate}
\end{Claim}


By an analogous argument, we have the following claim. 

\begin{Claim}\label{claim:neck_proof_2}
	For each $d\in D^1$, there is a sequence of coverings 
	$$P^\ast(x_d,t_d;r_d)\subseteq \bigcup_{i\in \tilde{N}^k_d}\mathcal{A}^i\cup \bigcup_{d\in \tilde{D}^k_d}P^\ast(x_d,t_d;r_d) \cup \bigcup_{e\in \tilde{E}^k_d}P^\ast(x_e,t_e;r_e)$$
	where $r_{d^\prime}\leq\frac{1}{2^{2k}}r_d$ for $d^\prime \in \tilde{D}^k$,
	\begin{enumerate}
		\item for every $i\in \tilde{N}^k_d$, $\mathcal{A}^i$ is an $(m,\epsilon,\gamma,\Lambda)$-annular region,
		\item for every $d^\prime \in \tilde{D}^k_d$, $P^\ast(x_{d^\prime},t_{d^\prime};r_{d^\prime})$ is a $(d)$-ball,
		\item for every $e\in \tilde{E}^k_d$, $P^\ast(x_e,t_e;r_e)$ is a $(e)$-ball,
		and we have the content estimate 
		$$ \sum_{i\in \tilde{N}^k_d}r_i^2+  \sum_{e\in \tilde{E}^k_d}r_e^2\leq C(Y,\gamma,\eta,\alpha,\Lambda)r_d^2.$$
	\end{enumerate}

\end{Claim}


Since $(M,(g_t)_{t\in [-10,10]})$ is smooth and compact, there is some $k_0\in \mathbb{N}$ such that the $k_0$-th element of any sequence constructed in Claims \ref{claim:neck_proof_1} and \ref{claim:neck_proof_2} must contain no $(c)$ or $(d)$-balls, respectively. Returning to our original cover \eqref{iden:neck_proof_cov_1}, we get after relabelling, a covering
$$ P^\ast(x_0,t_0;1)\subseteq \bigcup_{i\in N^1}\mathcal{A}^i\cup \bigcup_{e\in E^1}P^\ast(x_e,t_e;r_e),$$
where
$$\sum_{i\in N^1}r_i^2+\sum_{e\in E^1}r_e^2\leq C(Y,\eta,\gamma,\alpha,\delta,\Lambda).$$ By definition, for each $e\in E^1$ we have $$ \inf_{(x,t)\in P^\ast(x_e,t_e;\Lambda r_e,2r_e^2,2r_e^2)}\mathcal{N}_{(x,t)}(\delta r_e^2)\geq -Y+\delta,$$
and so we may repeat our construction on every $(e)$-ball in our cover (with constant $W=Y-\delta$) to produce a cover
$$P^\ast(x_e,t_e;r_e)\subseteq \bigcup_{i\in N^2_e}\mathcal{A}^i\cup \bigcup_{e^\prime \in E^2_e}P^\ast(x_{e^\prime },t_{e^\prime }; r_{e^\prime}).$$
Set $N^2=N^1\cup \bigcup_{e\in E^1}N^2_e$ and $E^2=\bigcup_{e\in E^1}E^2_e$, so that we now have a cover 
$$P^\ast(x_0,0;1)\subseteq \bigcup_{i\in N^2}\mathcal{A}^i \cup \bigcup _{e\in E^2}P^\ast(x_e,t_e;r_e)$$
where 
\begin{align*}
	\sum_{i\in N^2}r_i^2+\sum_{e\in E^2}r_e^2\leq C(Y,\gamma,\eta,\alpha,\delta,\Lambda)
\end{align*}
and 
\begin{align*}
 \inf_{(x,t)\in P^\ast(x_e,t_e;\Lambda r_e,2r^2_e,2r_e^2)}\mathcal{N}_{(x,t)}(\delta r_e^2)\geq -Y+2\delta
\end{align*}
for each $e\in E^2$. Iterating this process at most $Y\delta^{-1}$ times, we produce a cover
$$P^\ast(x_0,0;1)\subseteq \bigcup_{i\in N^{\lceil Y\delta^{-1} \rceil}}\mathcal{A}^i\cup \bigcup_{e\in E^{\lceil Y\delta^{-1} \rceil}} P^\ast(x_e,t_e;r_e)\subseteq P^\ast(x_0,0;2\Lambda,2,2),$$
with content bounds $$\sum_{i\in N^{\lceil Y\delta^{-1} \rceil}}r_i^2+\sum_{e\in E^{\lceil Y\delta^{-1} \rceil}}r_e^2\leq C(Y,\alpha,\eta,\gamma,\delta,\Lambda),$$
such that $$ \inf_{(x,t)\in P^\ast(x_e,t_e;\Lambda r_e,2 r_e^2,2r_e^2)}\mathcal{N}_{x,t}(\delta r_e^2)\geq -\sigma$$ 
for every $e\in E^{^{\lceil Y\delta^{-1} \rceil}}$. In particular, by \cite[Theorem 2.11]{Bam3}, if $\sigma \leq \overline{\sigma }(Y)$, we must have $r_{Rm}\vert_{P^\ast(x_e,t_e;r_e)}\geq c(Y,\delta)r_e$ for every $e\in E^{\lceil Y\delta^{-1}\rceil}$.
Let $\{\mathcal{B}^j\}_{j\in \hat{N}}$ denote the connected components of $\bigcup_{e\in E^{\lceil Y\delta^{-1}\rceil}}P^\ast(x_e,t_e;r_e)$. Clearly the covering we have produced satisfies the requirements of the theorem, and so the proof is complete.
\end{proof}

\section{Proof of the Main Theorems}\label{sec:main_thm}

Using Theorem \ref{thm:neck_decomp}, we prove Theorem \ref{thm:minkowski}.

\begin{proof}
   Fix $\Lambda\geq\underline{\Lambda}(Y)$, $\gamma\leq \overline{\gamma}(Y,\Lambda)$ as in Theorem \ref{thm:neck_decomp}, and take $\epsilon\leq\overline{\epsilon}(Y,\gamma,\Lambda)$ as in Lemma \ref{lem:annular}. We may apply Theorem \ref{thm:neck_decomp} to produce a cover  
    $$P^\ast(x_0,0;1)\subseteq \bigcup_{i\in N}\mathcal{A}^i\cup \bigcup_{j\in \hat{N}}\mathcal{B}^j\subseteq P^\ast(x_0,0;\Lambda,2,2)$$
    where
      \begin{enumerate}
        \item for each $i\in N$, $\mathcal{A}^i=P^\ast(x_i,t_i;2\Lambda r_i,2r_i,2r_i)\backslash \bigcup_{c\in\mathcal{C}^i}P^\ast(x_c,t_c;r_c)$ is an $(m,\epsilon,\gamma,\Lambda)$-annular region for some $m\leq e^{Y}$,
        \item for each $j\in\hat{N}$, we have  $\mathcal{B}^j=\bigcup_{e\in E^j}P^\ast(x_e,t_e;r_e)$  where for each $e\in E^j$ we have $r_{\operatorname{Rm}}\vert_{P^\ast(x_e,t_e;r_e)}\geq c(Y,\epsilon)r_e$,
    \end{enumerate}
     and we have the content estimate $$ \sum_{i\in N} r_i^2+ \sum_{j\in \hat{N}}\sum_{e\in E^j}r_e^2\leq C(Y,\epsilon).$$
      By construction, $\{(x,t)\in P^\ast(x_e,t_e;r_e)\,\vert\, r_{\operatorname{Rm}}(x,t)<\sigma \}$ is non-empty only if there $r_e<C(Y,\epsilon)\sigma$. Denote the collection of balls that satisfy this condition by $\hat{E}$. By Remark \ref{rem:annular_struc_1}, we have $$\operatorname{diam}_{P^\ast}(\{(x,s)\in \mathcal{A}^i\,\vert\, r_{\operatorname{Rm}}(x,s)<\sigma \}\cap (M\times \{t\})) \leq C(Y,\gamma,\Lambda)\sigma$$
 for each $i\in N$ and $t\in [\inf \mathfrak{t}(P^\ast(x_i,t_i;r_i)),\sup \mathfrak{t}(P^\ast(x_i,t_i;r_i))]$. Combining these facts with \cite[Theorem 9.8]{Bam1} we get that
   \begin{align*}   
   	| \{(x,t)\in P^\ast(x_0,0;1)\,\vert\, r_{\operatorname{Rm}}(x,t)<\sigma \}| & \leq \sum_{i\in N}\int_{t_i-r_i^2}^{t_i+r_i^2}|\{(x,s)\in \mathcal{A}^i\,\vert\, r_{\operatorname{Rm}}(x,s)<\sigma \}\cap (M\times \{t\})|_{g_t}\,dt\\&+\sum_{j\in \hat{E}}|P^\ast(x_e,t_e;r_e)|\\ &  \leq C(Y,\gamma,\epsilon)\sigma^4\left( \sum_{i\in N}r_i^2+\sum_{j\in \hat{N}}\sum_{e\in E^j}r_e^2\right)\\&\leq C(Y,\epsilon,\gamma)\sigma^4,
      \end{align*}
which completes the proof.
\end{proof}

Lastly, we prove Theorem \ref{thm:diffeo}.

\begin{proof}
     Choose $\Lambda\geq \underline{\Lambda}(Y)$ and $\gamma\leq \overline{\gamma}(Y,\Lambda)$ as in Theorem \ref{thm:neck_decomp}.  We may apply Theorem \ref{thm:neck_decomp} to produce a cover  
    $$P^\ast(x_0,t_0;1)\subseteq \bigcup_{i\in N}\mathcal{A}^i\cup \bigcup_{j\in \hat{N}}\mathcal{B}^j\subseteq P^\ast(x_0,t_0;2\Lambda,2,2)$$
    where
      \begin{enumerate}
        \item for each $i\in N$, $\mathcal{A}^i=P^\ast(x_i,t_i;2\Lambda r_i,2r_i,2r_i)\backslash \bigcup_{c\in\mathcal{C}_i}P^\ast(x_c,t_c;r_c)$ is an $(m,\epsilon,\gamma,\Lambda)$-annular region for some $m\leq e^{Y}$,
        \item for each $j\in\hat{N}$, we have  $\mathcal{B}^j=\bigcup_{e\in E^j}P^\ast(x_e,t_e;r_e)$  where for each $e\in E^j$ we have $r_{\operatorname{Rm}}\vert_{P^\ast(x_e,t_e;r_e)}\geq c(Y,\epsilon)r_e$,
    \end{enumerate}
     and we have the content estimate $$ \sum_{i\in N} r_i^2+ \sum_{j\in \hat{N}}\sum_{e\in E^j}r_e^2\leq C(Y,\epsilon).$$
     Set $E=\bigcup_{j\in \hat{N}}E^j$. We will need the following claim.
     \begin{Claim}\label{claim:main_thm_1}
         There exists $K\subseteq [-1,1]$ with $| K|<\epsilon $, such that for all $t\in [-1,1]\backslash K$, there exists $N_t\subseteq N$ and $E_t\subseteq E$ with $|E_t| ,|N_t|\leq C(Y,\epsilon)$, such that $$P^\ast_t(x_0,t_0;1)\subseteq \bigcup_{i\in N^1_t}\mathcal{A}^i_t\cup \bigcup_{e\in E_t}P^\ast_t(x_e,t_e;r_e)$$
         where $\mathcal{A}^i_t\coloneqq \mathcal{A}^i_t\cap (M\times \{t\})$.
     \end{Claim}
     \begin{proof}
    We  define $$f(x,t)=\sum_{i\in  N}r_i^{-4}\chi_{P^\ast(x_i,t_i;2\Lambda r_i)}+\sum_{e\in E}r_e^{-4}\chi_{P^\ast(x_e,t_e;2r_e)}$$
     where $\chi_{-}$ is the usual indicator function. By \cite[Theorem 9.8]{Bam1} we have 
     $$\int_{-1}^{1}\int_{P^\ast_t(x_0,0;2)} f(x,t)\,dg_tdt\leq C(Y)\left( \sum_{i\in N}r_i^2+\sum_{e\in E}r_e^2\right)\leq C(Y,\epsilon).$$
     In particular, we may choose a subset $K\subseteq [-1,1]$ with $|K|<\epsilon$ such that for all $t\in [-1,1]\backslash K$ we have $$\int_{P_t^\ast(x_0,0;2)}f(x,t)\,dg_t\leq C(Y,\epsilon).$$
     Moreover, arguing as in \cite[Proof of Lemma 6.1]{fluck2025epsilonregularitytheoremnoncollapsedricci}, we have
     $$\int_{P_t^\ast(x_0,0;2)}f(x,t)\,dg_t\geq c(Y)|\{i \in N\, \vert\, \mathcal{A}^i\cap (M\times \{t\})\neq \emptyset \}\cup \{j\in E\vert\, P^\ast(x_j,t_j;r_e)\cap (M\times \{t\})\neq \emptyset\}|,$$
     which completes the proof of the claim.
     \end{proof}

    For $t\in [-1,1]\backslash K$ as in Claim \ref{claim:main_thm_1}. Let  $\{\mathcal{B}^j_t\}_{ j\in \hat{N}_t}$ denote the connected components of $\bigcup_{e\in E_t}P_t^\ast(x_e,t_e;r_e)$. 
         
     \begin{Claim}\label{claim:final}
     For $\epsilon^\prime>0$ and $\epsilon\leq \overline{\epsilon}(\epsilon^\prime ,Y)$, the following is true. 	If  $t\in [-1,1]\backslash K $, $i\in N_t $, and $j\in \hat{N}_t$, there exists $\tilde{\mathcal{A}}^i_t\supseteq \mathcal{A}^i_t$ and $\tilde{\mathcal{B}}^j_t\supseteq \mathcal{B}^j_t$ such that the following is true. We have $$P^\ast_t(x_0,0;1)\subseteq \bigcup _{i\in N_t}\tilde{\mathcal{A}}^i_t \cup \bigcup_{j\in \hat{N}_t}\tilde{\mathcal{B}}^j_t\subseteq P^\ast_t(x_0,0;2\Lambda)$$
     	and \begin{enumerate}
     		\item $|N_t|+|\hat{N}_t|\leq C(Y,\epsilon)$,  \label{iden:neck_bubble_1}
     		\item for each $j\in \hat{N}_t$, we have $r_{Rm}\vert_{\tilde{\mathcal{B}}^j_t }\geq c(Y,\epsilon)\operatorname{diam}_{g_t}(\tilde{\mathcal{B}}^j_t)$,  \label{iden:neck_bubble_4}
         		\item for each $i\in N_t $, $\tilde{\mathcal{A}}^i_t$ is $\epsilon^\prime$-close in the $C^{\lfloor (\epsilon^\prime) ^{-1}\rfloor }$ topology to a standard metric annulus in a flat cone of the form $C(\mathbb{S}^3/\Gamma_i)$ for some $|\Gamma_i|\leq e^Y$, and $\tilde{\mathcal{A}}^i_t\cap \tilde{\mathcal{A}}^j_t \neq\emptyset$ if and only if $i=j$,\label{iden:neck_bubble_2}
     		\item for each $i\in N_t $ and $j\in \hat{N}_t$, the intersection $\tilde{\mathcal{A}}^i_t\cap \tilde{\mathcal{B}}^j_t$ is either empty or $\epsilon^\prime$-close in the $C^{\lfloor (\epsilon^\prime) ^{-1}\rfloor}$ topology to a standard metric annulus in a flat cone of the form $C(\mathbb{S}^3/\Gamma_i)$ for some $|\Gamma_i|\leq e^Y$.  \label{iden:neck_bubble_3}
     	\end{enumerate} 
     \end{Claim}
     
\begin{proof}
 For each $j\in N_t$, we let $\tilde{\mathcal{A}}^i_t=U_t$ be the open set constructed in Proposition \ref{prop:annular_time_slice}.  Note that \eqref{iden:neck_bubble_1} is immediate. If $\tilde{\mathcal{A}}^i_t\cap \tilde{\mathcal{A}}^j_t\neq\emptyset$, then necessarily $\Gamma_i=\Gamma_j$, and we may apply a standard gluing argument show that  $\tilde{\mathcal{A}}^i_t\cup \tilde{\mathcal{A}}^j_t$ is $\epsilon^\prime$-close in the $C^{\lfloor (\epsilon^\prime) ^{-1}\rfloor}$ topology to a standard metric annulus in the flat cone  $C(\mathbb{S}^3/\Gamma_i)$. Therefore, we may assume that $\tilde{\mathcal{A}}^i_t$ are pairwise disjoint, which verifies \eqref{iden:neck_bubble_2}.
 
We will now show \eqref{iden:neck_bubble_4}.  Fix $t\in [-1,1]\backslash K$ and consider some $j\in E_t$. 	Since $r_{Rm}$ is $1$-Lipschitz and $r_{Rm}\vert_{ P^\ast(x_{\hat{j}},t_{\hat{j}};r_{\hat{j}})  }\geq c(Y,\epsilon)r_{\hat{j}}$ for all  $\hat{j}\in E_t$, we may use \cite[Corollary 9.6]{Bam1} to deduce that if $\hat{j}\in E_t$ is such that $P^\ast(x_j,t_j;r_j)\cap P^\ast(x_{\hat{j}},t_{\hat{j}};r_{\hat{j}}) \neq \emptyset$ and $r_{\hat{j}}\leq c(Y,\epsilon)r_j$, then $ r_{Rm}\vert_{ P^\ast(x_{\hat{j}},t_{\hat{j}};r_{\hat{j}})  } \geq c(Y,\epsilon)r_{j}$. Therefore, any connected component $\mathcal{B}^j_t$ of $\bigcup_{e\in E_t}P^\ast(x_e,t_e;r_e)$ must satisfy $r_{Rm}\vert_{\mathcal{B}^j_t}\geq c(Y,\epsilon)\operatorname{diam}_{P^\ast}(\mathcal{B}^j_t)$. Combining \cite[Corollary 9.6]{Bam1} with some standard distance distortion estimates, we see that $\operatorname{diam}_{P^\ast}(\mathcal{B}^j_t)\geq c(Y,\epsilon)\operatorname{diam}_{g_t}(\mathcal{B}^j_t)$, which verifies \eqref{iden:neck_bubble_4}.
 
 Finally, we slightly modify the body regions $\mathcal{B}^j_t$ so that \eqref{iden:neck_bubble_3} is satisfied. Fix $j\in \hat{N}_t$. Given $i\in N_t$, if $\tilde{\mathcal{A}}^i_t\cap \mathcal{B}^j_t=\emptyset$  there is nothing to do,  so suppose that $\tilde{\mathcal{A}}^i_t\cap \mathcal{B}^j_t\neq \emptyset$. Suppose there exists $c^\prime > c$ such that $A_{cr_i,c^\prime r_i}\subseteq \varphi_{t}(\mathcal{B}^j_t)$ where $c=c(Y)$ and $\varphi_t$ are as in the proof of Proposition \ref{prop:annular_time_slice}. Let $c^\prime_0$ be the supremum over all  $c^\prime> c$ that satisfy this property.  Now remove $\tilde{\mathcal{A}}^i_t\cap \mathcal{B}^j_t\backslash \varphi_t^{-1}(A_{c^{\prime}_0r_i,10r_i})$ from $\mathcal{B}^j_t$. If such a $c^\prime$ does not exist, then there must exist $0<c^{\prime\prime}\leq 10$ such that $A_{c^{\prime\prime}r_i,10r_i}\supseteq \varphi_{t}(\mathcal{B}^j_t)$, and so  we remove $\tilde{\mathcal{A}}^i_t\cap \mathcal{B}^j_t\backslash \varphi_t^{-1}(A_{cr_i,c^{\prime\prime}_0r_i})$ where $c^{\prime\prime}_0>0$ is the infimum over all  such $0<c^{\prime\prime}\leq 10$.  We let $\tilde{\mathcal{B}}^j_t$ denote the region produced after applying this procedure for every $i\in N_t$. It is easy to see that $\tilde{\mathcal{A}}^i_t$  and $\tilde{\mathcal{B}}^j_t$ satisfy \eqref{iden:neck_bubble_3}, which completes the proof of the claim.
\end{proof}
       
Note that given $\epsilon^\prime>0$, we may assume without loss of generailty that $\epsilon\leq \overline{\epsilon}(\epsilon^\prime ,Y)$ as in Claim \ref{claim:final}. As in \cite[Section 8.4]{regEinstein}, we may choose $\epsilon^\prime\leq \overline{\epsilon^\prime}(Y) $ so that the attaching maps in Claim \ref{claim:final} are isotopic to an isometry, and so the proof follows in an identical fashion.

\end{proof}

\printbibliography

\end{document}